\documentclass[12pt]{article}

\usepackage{amsmath,amssymb,amscd,graphicx,epsfig}

\everymath{\displaystyle}

\newtheorem{teo}{Theorem}[section]
\newtheorem{lem}[teo]{Lemma}
\newtheorem{cor}[teo]{Corollary}

\newtheorem{prop}[teo]{Proposition}
\newtheorem{lem-defi}[teo]{Lemma-Definition}
\newtheorem{defi}[teo]{Definition}

\newtheorem{remark}[teo]{Remark}

\newcommand{\mr}{\mathbb{R}}
\newcommand{\mc}{\mathbb{C}}
\newcommand{\mz}{\mathbb{Z}}

\newcommand{\mn}{\mathbb{N}}

\newcommand{\ma}{\mathbb{A}}

\newtheorem{theo}{Theorem}[section]

\newtheorem{lemma}[theo]{Lemma}

\newcommand{\Aa}{{\mathcal A}}
\newcommand{\Bb}{{\mathcal B}}
\newcommand{\Cc}{{\mathcal C}}

\newcommand{\Kk}{{\mathcal K}}
\newcommand{\Ll}{{\mathcal L}}

\newcommand{\CM}{{\mathbb C}}
\newcommand{\NM}{{\mathbb N}}

\newcommand{\RM}{{\mathbb R}}

\newcommand{\ZM}{{\mathbb Z}}

\newcommand{\tr}{\mbox{\rm Tr}}                   
\newcommand{\TV}{{\mathcal T}}                    
\newcommand{\TVm}{{\mathcal T}_{\mu}}             
\newcommand{\Cs}{$C^{\ast}$-algebra }             
\newcommand{\Css}{$C^{\ast}$-algebras }           

\def\F{{\cal F}}

\def\A{{{\cal A}}}
\def\C{{{\cal C}}}

\def\H{{{\cal H}}}

\def\L{{{\cal L}}}

\def\P{{{\cal P}}}

\def\M{{{\cal M}}}

\def\F{{\cal F}}

\def\A{{{\cal A}}}
\def\C{{{\cal C}}}

\def\H{{{\cal H}}}

\def\L{{{\cal L}}}

\def\P{{{\cal P}}}

\title{Spaces of Tilings, Finite Telescopic Approximations  and
Gap-Labelling }

\author{\normalsize Jean Bellissard$^{1,2}$, Riccardo Benedetti$^{3}$, 
Jean-Marc
Gambaudo$^{4}$\\ {\small $^1$ Universit\'e Paul-Sabatier, 118
route de Narbonne, 31062 Toulouse, France}\\ {\small $^2$ Institut
Universitaire de France}\\ {\small $^3$ Dipartimento di
Matematica, via F. Buonarroti 2, 56127 Pisa, Italy }
\\ {\small $^4$ Laboratoire de Topologie, U.M.R. 5584 du CNRS,} \\ {\small 
Universit\'e de Bourgogne,
B.P. 47870- 21078 Dijon Cedex France} \vspace{.2cm}}

\date{ }

\begin{document}

\maketitle

\begin{abstract} For a large class of tilings of $\RM^d$, including the Penrose 
tiling in  dimension 2 as well as the icosahedral ones in  dimension 3, the continuous 
hull $\Omega_T$ of such a tiling $T$  inherits a minimal $\RM^d$-lamination 
structure with flat leaves 
and a  transversal $\Gamma_T$ which is a Cantor set. In this case, we show that 
the continuous hull can be seen as the  projective limit of a suitable sequence of 
branched,  oriented and flat compact $d$-manifolds. Truncated sequences
furnish better and better finite approximations of the asymptotic dynamical
system and the algebraic topological features related to this sequence reflect 
the dynamical properties of the $\RM^d$-action on the continuous hull. 
In particular the set of positive invariant measures of this action turns to be a convex 
cone, canonically associated with the orientation,
in the projective limit of the $d^{th}$-homology groups of the branched manifolds. 
As an application of this construction we prove a {\it gap-labelling
theorem}: 

Consider  the $C^*$-algebra $\Aa_T$ of $\Omega_T$, and the group 
$K_0(\mathcal{A}_T)$, then 
for every  finite $\RM^d$-invariant measure $\mu$ on $\Omega_{T}$,   the image 
of the group 
$K_0(\mathcal{A}_T)$ by the 
$\mu$-trace satisfies:
$$\TVm (K_0(\mathcal{A}_T))   \,=\, 
   \int_{\Gamma_T} d\mu^t \; \Cc (\Gamma_T, \ZM),$$ where $\mu^t$ is the 
transverse invariant measure on 
$\Gamma_T$ induced by $\mu$ and $\Cc (\Gamma_T, \ZM)$ is the set of continuous 
functions on $\Gamma_T$ with 
integer values.     
\end{abstract}

\section {Introduction}
\label{sec-intro}

\noindent It has been argued \cite{JB86,JB93} that the mathematical description 
of aperiodic solids can be made through the construction of the so-called {\em 
Noncommutative Brillouin Zone} (NCBZ) to replace the Bloch theory used for 
periodic crystals. In a recent work \cite{BHZ00}, the construction of this 
noncommutative manifold was performed from the equilibrium positions of the 
atoms of the solid, seen as a {\em uniformly discrete} subset of the physical 
space $\RM^d$. A suitable compactification of the family of $\RM^d$-translates 
of this set leads to the notion of {\em Hull}, which is a compact space $\Omega$ 
endowed with an action of $\RM^d$ by homeomorphisms. The \Cs of continuous 
functions on the NCBZ is then given by the {\em crossed product} \cite{Ped}
$\Cc (\Omega)\rtimes \RM^d$ of the algebra of the Hull by the $\RM^d$-action. 

The Hull can also
be seen as a {\em lamination} \cite{Ghys} or a {\em foliated space} \cite{MoSc}, 
namely a foliation with non smooth 
transverse structure. On the other hand \cite{tilings}, the 
construction of the Voronoi cells from the point set of atomic positions, leads 
to a {\em tiling} of $\RM^d$ by polyhedra touching face to face, from which the point set can be 
recovered by a dual construction. Hence, the construction of the Hull can 
equivalently be performed from three complementary point of views: as a 
dynamical system, as a lamination or foliated space, as a tiling. 
This latter point of view permits to select constraints using the tiling 
language more easily than using the language of point sets. The tilings that 
have mostly attracted the attention of experts are those with a finite number 
of patches with a bounded size modulo translations (the so called {\em finite pattern condition} or FPC) 
and satisfying {\em repetitivity}. 
This imposes severe constraints on the point set of 
atomic positions. In particular, repetitivity implies the minimality of the
$\RM^d$-action, while the FPC implies that the lamination is transversally a Cantor set. 
It is fair to say however, that some tilings of interest have 
a finite number of tiles only if we identify tiles under both translation and 
rotation \cite{pinwheel}: these are not considered in the present paper.

One of the main consequence of constructing the NCBZ is the so-called {\em Gap 
Labelling Theorem} \cite{JB82,JB86,JB93}. If $\Aa$ is the \Cs of continuous 
functions on the NCBZ, and if $H=H^{\ast}$ is a selfadjoint element of $\Aa$, 
any spectral gap of $H$ can be associated with a canonical element of the 
$K$-group $K_0(\Aa)$ \cite{Black86}. Moreover, the $\RM^d$-invariant ergodic 
probability measures $\mu$ on the Hull are in one-to-one correspondence with the 
extremal {\em traces} on $\Aa$. Each such trace $\tau_{\mu}$ canonically defines 
a group homomorphism $\tau_{\mu,\ast}$ from $K_0(\Aa)$ into $\RM$ the image of 
which is called the set of {\em gap labels}. Since $\Aa$ is separable it follows 
that the set of gap labels is a countable subgroup of $\RM$. It has been 
conjectured \cite{BHZ00} that, if the Hull is completely disconnected 
transversally to the $\RM^d$ action, then the set of gap labels is nothing but 
the $\ZM$-module generated by the integrals $\int_{\Omega} d\mu\, f$ whenever $f$ runs 
through the set $\Cc (\Omega, \ZM)$ of integer valued continuous functions on 
the Hull. This conjecture was proved for the first time for $d=1$ in \cite{JB93} 
using the Pimsner-Voiculescu exact sequence \cite{PV80}. The same method was 
used in \cite{VE94} to extend the result for $d=2$. In \cite{84}, the result for 
$d=2$ was reestablished by using the Kasparov spectral sequence. 
In the case for which the Hull is given by a $\ZM^d$-action on a Cantor set $X$, 
Hunton and Forrest \cite{FH99}, using the technic of spectral sequences,  have 
proved that the $K$-group is isomorphic to the group cohomology of $\ZM^d$ with 
values in the group $\Cc(X,\ZM)$. While this result does not lead to the computation of 
the set of gap labels in general, it permits to compute the $K$-group 
in many practical situations that occur in physics \cite{KHF,KG} as well as the 
set of gap labels for $d=3$ \cite{BKS01}. The proof of this conjecture 
for tilings that satisfy repetitivity and FPC in an arbitrary dimension is one of the main results of this paper. 

\begin{remark} 
{\rm At the time this paper was written, 
Benameur and Oyono \cite{BenOy} on the one hand, Kaminker and Putnam 
\cite{KamPut} on the other, announced that they have proved the gap labelling 
conjecture for a $\ZM^d$-action on a Cantor set for any $d$. Our result goes
beyond such situation.}
\end{remark}

\vspace{.2cm}

In the present work, the construction of the Hull for tilings with a finite 
number of tiles modulo translations, is given in term of an inverse limit or 
{\em expanding flattening sequence} (EFS) of {\em branched oriented flat 
manifolds with dimension $d$} (BOF-$d$-manifolds), called the approximants. The existence of a flat metric on each
BOF-manifold leads to the notion of parallel transport by constant vector fields and to an
$\RM^d$-action on the inverse limit. Consequently, the {\em $\RM^d$-action on the Hull is topologically 
conjugate to the $\RM^d$-action on inverse limit of the approximants}.
This construction was inspired by previous works on dynamical 
systems which can be seen as tilings in dimension 1 
\cite{Versik,Putnam,G-M}.  

\noindent One consequence of this construction is the identification of the space of
positive invariant measures on the Hull with the positive cone of the inverse limit of the top-homology 
groups of the approximants. This cone is canonically 
defined thanks to the orientation of the approximants. Hence, the set of positive invariant
measures can be interpreted from various equivalent point of views: as the positive traces
on $\Aa$, as the positive cone of the $d$-homology of the Hull, or, as will be seen later,
as {\em Ruelle-Sullivan currents} on the lamination.

\noindent The \Cs of the noncommutative Brillouin zone can be constructed from
the inductive limit of the space of 
continuous functions on the approximants. 
Consequently, thanks to the Thom-Connes theorem, the $K_0$-group of the
NCBZ can be seen as the inductive limit of the $K_0$-group of the approximants.
For a BOF-$d$-manifold, the $K$-theory can be computed from a spectral
sequence very similar to the one 
constructed for $\ZM^d$-action on a Cantor set \cite{FH99}. Moreover, the set of 
gap labels can also be obtained from the so-called {\em Ruelle-Sullivan current} 
of the lamination, which in turn is shown in this paper to be given by a 
compatible family of $d$-cycles defined on each BOF-manifold of the EFS
defining the lamination. The explicit construction of such family is given here 
and permits to prove the gap labelling theorem for the class of tilings under 
investigation.

\subsection {Tilings}
\label{ssec-tilings}

\noindent This subsection is devoted to few important properties of tilings in 
 the usual $d$-dimensional Euclidean space
$\RM^d$. For any 
$x$ in $\RM^d$ and any $\epsilon >0$, $B_\epsilon(x)$ denotes the closed
$d$-ball in $\RM^d$ with center $x$ and radius $\epsilon$ and 
$\Vert \cdot \Vert$ the Euclidean norm in $\RM^d$. The vector space $\RM^d$ is
{\it oriented} by stipulating that the standard basis is positive.  A
{\it prototile} is a compact closed connected  subset of ${\RM^d}$ which is the closure of its interior{\footnote{It is often assumed in the litterature that prototiles are homeomorphic to closed ball. The  genralization we impose here is motivated by further constructions, but the standard results we will recall, work as well in our situation.} A prototile determines, up to  direct isometries  of ${\RM^d}$, a {\it tile $i$-type} and up to translations  of ${\RM^d}$ only, a {\it tile $t$-type}. Once the kind of tile type $i$ or $t$ is fixed, a prototile will often be identified with its tile type. Let  
$X$ be a {\it countable} collection of prototiles of pairwise distinct types (of a given kind $i$ or $t$).  A tiling 
$T$ with tile types in $X$, consists of a countable collection of subsets of $\RM^d$, $\{t_1, 
t_2,\dots, t_n,\dots\}$, called {\it tiles}, such that:

\begin{itemize}
  \item each tile belongs to a tile type in  $X$;
  \item the union of the tiles $\cup_{i\geq 1}t_i$ covers the whole
  space $\RM^d$;
  \item the interiors of the tiles are pairwise disjoint.
\end{itemize}

\noindent   A tiling $T$ is {\it polyhedral} if all its tiles are   polyhedra in $\RM^d$. 

\noindent An example of polyhedral
tiling is the Voronoi tiling of a Delone set \cite{Lag}. A countable subset $\Ll$ of $\RM^d$
is  {\em uniformly discrete} if there is $r>0$ such that every open ball of radius $r$
meets $\Ll$ on at most one point. $\Ll$ is {\em relatively dense} if there is $R>0$
such that every open ball of radius $R$ meets $\Ll$ on at least one point. 
$\Ll$ is a {\em Delone set} if it is both uniformly discrete and relatively dense.
If $x\in \Ll$, the {\em Voronoi cell} $V_x$ is the open set of points in $\RM^d$
that are closer to $x$ than any other points in $\Ll$. If $\Ll$ is a 
Delone set, each Voronoi cell contains an open ball of radius $r$ and is contained into 
an open ball of radius $R$. Moreover, the family $\{ \overline{V_x}\,; \; x\in \Ll \}$
defines a polyhedral tiling of $\RM^d$, in which the tiles are convex and meet face-to-face. 
Conversely, let $X$ be a set of prototiles with pairwise distinct tile types, such
that each prototile contains a ball of radius $r$ and is contained in a ball of radius $R$. 
Define in each prototile $p$ of $X$ a point $x_p$, namely {\em $X$ is punctured}. 
Then if $T$ is a tiling  with type tiles in $X$, the family of points 
$x_t$ associated with each tile $t$ by identification of $t$ with a prototile, is a Delone set. 
However, the corresponding Voronoi tiling does not usually coincide with $T$. 

\noindent Given a countable set $X$ of prototiles, it is not known in general
whether the set of tilings with type tiles in $X$ is empty or not.
However, many examples of finite $X$'s are known to give a rich set of tilings. For example, tiles of 
the Penrose tiling or the Voronoi tiling of a quasi crystal \cite{tilings}, 
belong to finitely many  $t$-types in $X$, while tiles of the pinwheel tiling \cite{pinwheel} belong  only to  
two tile $i$-types but to  an infinite number of tile $t$-types. 

 {\sl In this paper we focus our attention on  the following type of tiling spaces:
\begin{itemize}
\item we work with tile $t$-type and assume that the set of prototiles $X$ is finite;
\item we consider the tiling space $T(X)$ made with tilings for which the tile $t$-type of each tile is in $X$;
\item furthermore, it will be always assumed that $T(X)\neq \emptyset$.
\end{itemize}
}

  The group 
$(\RM^d,+)$ acts on $T(X)$ by translations:

$$\omega : \RM^d \times T(X)\to T(X)\ \ \ \ \omega (x,T)=T+x.$$ 

\noindent The tiling space $T(X)$ is 
endowed with a distance (hence with the induced topology) defined as follows: 
let $A$ denote the set of $\epsilon \in ]0,1[$ such that there exist $x$ and 
$x'$ in  $B_\epsilon(0)$ such that $(T+x)\cap B_{1/\epsilon}(0) = (T'+x')\cap 
B_{1/\epsilon}(0)$, then:

$$ \delta (T,T')= \inf A \ \ \ \ {\rm if} \ \ A\neq \emptyset$$

$$\delta (T,T')=1  \ \ \ \ {\rm if} \ \ A = \emptyset \ .$$

\noindent Hence the diameter of  $T(X)$ is bounded by $1$ and the $\RM^d$-action  
on $T(X)$ is continuous. 

\noindent For each $T\in T(X)$, $o(T)= T+\RM^d$ denotes its {\it 
orbit}. The distance $\delta$ restricts to any orbit $o(T)$ and the induced topology 
is finer that the one induced by the Euclidean distance in $\RM^d$; this 
topology  reflects the way two translated tilings look the same at a local 
level. The  {\it continuous Hull} of $T,\,$ $\Omega_T, $ is the closure of $o(T)$ in $T(X)$.  It is invariant for the $\RM^d$-action on $T(X)$.

\begin{defi}
\label{fpc}
{\rm A {\it patch} of $T\in T(X)$ is a finite sub-collection of tiles of $T$.
The tiling $T$ satisfies the {\it finite pattern condition} FPC if for any 
$s>0$, there are {\it up to translation}, only finitely many patches with
diameter smaller than $s$. }
\end{defi}
We refer to  \cite{Lag,BHZ00} to understand, what means the FPC from the point set point of view.

\begin{remark}
{\rm 
Actually, the FPC is automatically satisfied in the case of a polyhedral tiling of $\RM^d$, in which the prototiles are in finite number and  meet face-to-face. However in the general case it is easy to construct tilings with a finite number of prototiles which do not satisfy the FPC.}
\end{remark}

\noindent This last definition is motivated by the following standard result
(see the body and the references of \cite{KP00}).

\begin{prop}
\label{omegacomp} 
If a tiling $T\in T(X)$ satisfies the FPC, then $\Omega_T$ is compact so it 
coincides with the metric completion of $(o(T), \delta \upharpoonright_{o(T)})$.
\end{prop}

\begin{remark}
{\rm If the FPC holds, the metric on $T(X)$ induces the weak$^{\ast}$ topology
defined in \cite{BHZ00}. In other cases, while the tiling distance might
lead to non compact continuous Hull, the weak$^{\ast}$ topology always lead to a compact Hull.  }
\end{remark}

\noindent Notice  that if a tiling $T$ satisfies the finite pattern condition, 
then all tilings in $\Omega_T$ satisfy also the finite pattern condition.

\begin{defi}
\label{transversale} 
{\rm Let $X_p$ be a punctured copy of $X$, that is each prototile $p\in X$ is given with a marked point $x_p$ in its interior. Consequently  each tile $t$ in a tiling $T\in T(X)$ admits a distinguished point $x_t$.
The set $\Sigma_{T_p}$ of tilings $T' \in \Omega_T$ such that one of the $x_t$
coincides with the origin $0\in\RM^d$  is called the {\em canonical transversal}
associated with $X$. }
\end{defi}

\noindent The following result is also standard \cite{KP00} and admits an
equivalent for Delone sets of finite types \cite{BHZ00}.

\begin{prop}
\label{cantortransve} 
If a tiling $T\in T(X)$ satisfies the FPC,then for every punctured $X_p$,  the canonical
transversal $\Sigma_{T_p}$ is compact and completely disconnected.
\end{prop}

\begin{defi}
\label{ripete} 
{\rm A tiling $T\in T(X)$ satisfies the {\it repetitivity}
condition if for any patch in $T$ there exists $s>0$, such that
for every $x$ in $\RM^d$, there exists a translate of this patch which
is in $T$ and in the ball $B_s(x)$.}
\end{defi}

\noindent Again this last definition is motivated by the following
proposition:

\begin{prop} If a tiling $T\in T(X)$ satisfies both finite pattern and
repetitivity conditions then the orbit $o(T')$ of any $T'\in
\Omega_T$ is dense in $\Omega_T$. In other words the dynamical
system $(\Omega_T, \omega)$ is {\rm minimal}.
\end{prop}

\begin{defi}
\label{aperiodic}
{\rm A tiling $T\in T(X)$ is {\it aperiodic} if
there exists
no $x\neq 0$ in $\RM^d$ such that $T+x =T$; it is {\it strongly
aperiodic} if all tilings in $\Omega_T$ are aperiodic.}
\end{defi}

\noindent We have the following classical result (see \cite{KP00}):

\begin{prop}\label{rectif}
If an aperiodic tiling $T\in T(X)$ satisfies both finite
pattern and repetitivity conditions then it is strongly aperiodic. In this case,
 any canonical
transversal $\Sigma_{T_p}$ is a Cantor set.

\end{prop}

\begin{defi}\label{perfect}{\rm The tilings $T\in T(X)$  which are
aperiodic and
satisfy both finite pattern and repetitivity conditions are said
{\it perfect} tilings.}
\end{defi}

\noindent As the main object is the dynamical system
$(\Omega_T,\omega)$, this suggests the following equivalence relation
on tilings of $\RM^d$.

\begin{defi}
\label{omegaequiv}
{\rm Two tilings of $\RM^d$, $T$ and $T'$,
are {\it $\Omega$-equivalent} if there exists a
homeomorphism $\phi: \Omega_T\to \Omega_{T'}$ which conjugates the two
$\RM^d$-actions.}
\end{defi}

\subsection{Main results}
\label{ssec-content}

\noindent Let us describe now the content of the paper.

\vspace{.1cm}

\noindent $\bullet$ In section
\ref{sec-TvsL} some basic general notions about {\it
laminations} are given and, for a perfect tiling $T$ of $\RM^d$, it is proved that 
$\Omega_T$ can be equipped with a natural lamination structure
whose leaves correspond to the $\RM^d$-orbits and which is
transversely a Cantor set. Those laminations arising from perfect tilings are 
called {\em tilable} and are characterized.  As an
application we prove the following result:
\begin{teo}Any perfect tiling is
$\Omega$-equivalent to some perfect  tiling whose prototiles are $d$-rectangles. 
\end{teo} 

\vspace{.1cm}

\noindent $\bullet$ In section \ref{sec-BOFS} the category of branched
oriented and flat (BOF) $d$-manifolds is defined. Local models, 
morphisms and in particular {\it BOF-submersions}, homology and
cohomology are defined and studied.

\vspace{.1cm}

\noindent $\bullet$  A particular class of projective limits of BOF-$d$-manifolds
is defined in section \ref{ssec-EFS}. These projective limits are defined through
BOF-submersions from one BOF-$d$-manifold to the previous one satisfying a
{\it flattening property} introduced by  Williams in \cite{Williams}. The sequences
that form the projective limits in this class are called {\it expanding 
flattening sequences} (EFS). Here again the homological and cohomological
properties of these projective limits are described. 

\vspace{.1cm}

\noindent $\bullet$  In section \ref{ssec-JM} the connection between expanding 
flattening sequences and tilings is analyzed and leads to the following result:

\begin{teo}
Let $T$ be a perfect tiling. Then its continuous Hull,
seen as a dynamical system $(\Omega_T, \omega)$, 
is topologically conjugate to the inverse limit of an expanding flattening sequence
of branched oriented flat manifolds, endowed with the action of $\RM^d$ by parallel transport 
under constant vector fields. 
\end{teo}

\noindent On  one hand,  an EFS defines a 
compact $\RM^d$-invariant subset of $T(X)$, for some finite set $X$ of $d$-rectangles
prototiles.  Further conditions on the EFS ensure that this inverse limit is just one 
(perfect) $\Omega_T$. On the other hand, for every perfect  tiling $T$, 
the dynamical system $(\Omega_T,\omega)$ can be realized as the inverse limit 
of a suitable EFS.  Furthermore, not only the topological space $\Omega_T$ can be 
obtained by this inverse limit procedure but also most of the relevant structures on $\Omega_T$.
Hence, the $\RM^d$-action, the topological $K$-theory groups $K^i(\Omega_T)$ 
and so on, can be obtained as inverse (or possibly direct) limits of suitable 
``pre-structures'' on each of the BOF-$d$-manifolds forming the EFS. The ``finite
telescopic approximations'' of $(\Omega_T,\omega)$ mentioned in
the title are the {\it truncated} finite subsequences of the
EFS, decorated with their truncated sequences of pre-structures.
These truncated sequences can be encoded by means of finite
combinatorial patterns and provide a better and better approximation of the
asymptotic dynamical system $(\Omega_T,\omega)$ as the
length of the truncated sequences grows up. The inverse limit
approach to perfect {\it substitution tilings} \cite{KP00}
 can be rephrased in terms of EFS and regarded as a
special case of the previous construction. For substitution tilings the
BOF-$d$-manifolds of the EFS are all equal each to the other, up to an  overall
rescaling. 

\vspace{.1cm}

\noindent $\bullet$  A particular attention is paid to the ergodic theory of 
tiling spaces in section \ref{invariantmes}. 
In the general case of a lamination  with $d$-dimensional leaves,  
it is known that a transverse invariant measure $\mu^t$ can be seen as 
a Ruelle-Sullivan current i.e. a linear form on the 
vector space of $d$-forms on the lamination which satisfies some extra 
properties (see section \ref{invariantmes}). In the particular case of a perfect tiling space 
$(\Omega_T,\omega)$ the correspondence between  the transverse invariant measures 
of the lamination, the invariant measures of 
the $\RM^d$-action and the invariant measure on a transversal under the action 
of the holonomy groupoid  of the associated lamination are also recalled. 
Considering $\Omega_T$ as a projective limit of an EFS of 
BOF-$d$-manifolds $B_i$, the data of a finite $\RM^d$-invariant measure $\mu$  
on $(\Omega_T,\omega)$ correspond to a {\it positive weights system}  $\mu_i$
on the $d$-cells of a natural cell decomposition of $B_i$,
satisfying a Kirchoff-like law. An elementary, but important fact is that the
orientation, which is part of the BOF structure, allows 
to interpret these weights systems also as a positive cone in the space of  
$d$-cycles on $B_i$. So it is
possible to contract these cycles against any real $d$-cohomology
class. This passes to the limit and the set of finite invariant measures of 
$(\Omega_T,\omega)$ 
corresponds to the convex set of the projective limit of the positive cones of 
the $d^{th}$ homology groups of the BOF-$d$-manifolds $B_i$.  On the other hand, 
 the $d$-cohomology group on $\Omega_T$ can be seen as the inductive limit of 
the $d$-cohomology groups of the BOF-$d$-manifolds $B_i$, each real
$d$-cohomology class $\alpha$ on $\Omega_T$ is the direct limit of classes 
$\alpha_i$ on $B_i$. 
It follows that for every finite invariant measure $\mu$, the following pairing can be defined:

$$<\mu|\alpha>\, =\,<\mu_i\vert \alpha_i>\,\in \mr,$$

\noindent for $i$ big enough. The special case of integer valued cohomology classes 
is studied in more details in view of the gap-labelling theorem.

\vspace{.1cm}

\noindent $\bullet$  Section \ref{sec-cstar} is devoted to the non-commutative 
aspects of tiling spaces 
theory with a special emphasis on a gap labelling theorem. The fundamental
 non-commutative geometric structure associated with the dynamical
system $(\Omega_T,\omega)$ is a suitable $C^{\star}$-algebra
$\mathcal{A}_T$.  The topology of this non commutative manifold is given by
the group $K(\mathcal{A}_T) = K_0(\mathcal{A}_T) \oplus K_1(\mathcal{A}_T) $ \cite{Black86} (by Bott periodicity theorem $K_i \cong K_{i+2}$).
The Thom-Connes theorem \cite{Black86} shows that 
$K_i(\mathcal{A}_T)$ is isomorphic to the topological $K$-group $K^{i+d}(\Omega_T)$.
On the other hand, with each finite $\RM^d$-invariant 
measure $\mu$ on $\Omega_T$  is associated a canonical {\em trace} on $\Aa_T$ which induces 
a group homomorphism $\TVm: K_0(\mathcal{A}_T)\cong 
K^d(\Omega_T)\to \RM$. The study  of the image of $\TVm$ is relevant in
the framework of the so called gap labelling \cite{JB82,JB86,JB93,KP00,BHZ00,BKS01}. 
The isomorphism $ K_0(\mathcal{A}_T)\cong
K^d(\Omega_T)$ allows to factorize the pairing between invariant measures and 
elements in $K_0(\mathcal{A}_T)$ has a pairing between a cycle (the invariant measure) and an 
integer valued cocycle (see section \ref{sec-cstar} for further precision), leading to:

\begin{teo} {\rm (Gap-labelling theorem)}
\label{gaplab} 
Let $(\Omega_T,\omega)$ be the dynamical system associated with a perfect tiling 
$T$.  Let $\mu$ be a finite invariant measure, $\Gamma_T$ be a transversal and $\mu^t$ 
the induced transverse invariant measure on $\Gamma_T$, then:

$$\TVm (K_0(\mathcal{A}_T))   \,=\, 
   \int_{\Gamma_T} d\mu^t \; \Cc (\Gamma_T, \ZM),$$ 

\noindent where  $\Cc (\Gamma_T, \ZM)$ is the set of continuous functions on 
$\Gamma_T$ with integer values.  
\end{teo}

\begin{remark}
{\rm
\noindent 1)- The existence of a underlying lamination structure on
$\Omega_T$ is implicit in the literature about tiling dynamics and, for
instance, it is explicitely described in \cite{Ghys,MoSc}.  The idea of getting
finite approximations  by means of suitable
branched manifolds has, as a natural background, the encoding of
measured laminations embedded into 2-dimensional (3-dimensional)
manifolds by means of embedded ``train-tracks'' (``branched
surfaces'') equipped with positive weight systems \cite{Williams,T}. 
In the
realm of tiling dynamics, it has been already
observed in \cite{Versik} and used in \cite{Putnam,G-M},
that the dynamics of a minimal map on a Cantor set can be described as
an inverse limit of oriented graphs. On the other hand, the dynamical
systems induced by one dimensional perfect tilings of $\RM$  exactly
corresponds to minimal dynamics on a Cantor set. Despite the fact that
in dimension greater than $1$, tilings are much more
complicated objects, it turns out that this approach can be
generalized by using the convenient category of branched oriented flat
manifolds.

\vspace{.1cm}

\noindent 2)- In the present paper only tilings with a finite number of $t$-types are considered so that $\Omega_T$ is in a natural way a dynamical system with an action of $\RM^d$ by translations. One could ask whether our constructions could be generalized to the case of tilings with a finite number of $i$-types and with hulls carrying an action of the full group of direct
isometries of $\RM^d$ could also be considered. We hope to face this problem in a
forthcoming  work.
}
\end{remark}

\noindent{\bf Acknowledgments:} It is a pleasure for J.-M. Gambaudo to
thank I. Putnam for very helpful (e-mail) conversations, and in particular for having pointed out a mistake in a first version of the paper. He also thank the
Department of Mathematics of the University of Pisa, where a part of this
work has been done, for its warm hospitality. J. Bellissard is indebted to 
several colleagues, A. Connes, T. Fack, J. Hunton, J. Kellendonk, A. Legrand, 
who helped him becoming more familiar with the technicalities involved in the 
gap labelling theorem. He thanks M. Benameur and I. Putnam for letting him know
their recent works \cite{BenOy,KamPut} prior to publication. He also wants to 
thank the MSRI  (Berkeley) and Department of Mathematics of the University of 
California at Berkeley for providing him help while this work was done during the 
year 2000-2001. 

\section {Tilings versus Laminations}
\label{sec-TvsL}

\noindent Let $M$ be a compact metric space and assume there exist a cover
of $M$ by open sets $U_i$ and homeomorphisms called {\it charts}
$h_i:U_i\to V_i\times T_i$ where $V_i$ is an open set in $\RM^d$
and $T_i$ is some topological space. These open sets and
homeomorphisms define an atlas of a {\it $(d)$-lamination}
structure with $d$-dimensional leaves on $M$, if the {\it
transition maps} $h_{i,j}\, =\, h_j\circ h_i^{-1}$ read on their
domains of definitions: $$h_{i, j}(x, t)\,=\, (f_{i,j}(x, t),\,
\gamma_{i, j}(t)),$$ where $f_{i,j}$ and $\gamma_{i, j}$ are
continuous in the $t$ variable and $f_{i, j}$ is
smooth in the $x$ variable.  Two atlas are {\it
equivalent} if their union is again an atlas. A {\it $(d)$-lamination}
is the data of a metric compact space $M$ together
with an equivalence class of atlas $\mathcal{L}$.

We call {\it slice} of a lamination a subset of the form
$h_i^{-1}(V_i\times \{t\})$. Notice  that from the very definition of a lamination, a slice associated with some chart $h_i$ intersects at most one slice associated with another chart $h_j$. The {\it leaves} of the lamination
are the smallest connected subsets that contain all the slices
they intersect. Leaves of a lamination inherit a d-manifold
structure, thus at any point in the lamination, it is possible to
define the tangent space to the leaf passing through this point. A
lamination is {\it orientable} if there exists in $\mathcal{L}$ an
atlas made of charts $h_i:U_i\to V_i\times T_i$ and orientations
associated with each $V_i$ preserved by the restrictions $f_{i, j}$
of the transition maps to the leaves. $\mathcal{L}$ is {\it
oriented} if one fixes one global orientation. 
\begin{defi}
\label{def-trans}
 {\rm Given a lamination
$(M,\mathcal{L})$, a {\it transversal} of $\mathcal{L}$ is a
compact subset $\Gamma$ of $M$ such that, for any leaf $L$ of
$\mathcal{L}$, $\Gamma \cap L$ is non empty and is a discrete
subset with respect to the $d$-manifold topology of the leaf $L$.}
\end{defi}

\noindent The laminations which are related to tilings are in fact of a very
special kind.

\begin{defi}
\label{flatlam} 
{\rm An oriented $(d)$-lamination$(M, \mathcal{L})$
is {\it flat} if the following properties are satisfied:

\noindent 1) there exists in $\mathcal{L}$ a maximal atlas made of charts
$h_i:U_i\to V_i\times T_i$ such that  the transition maps  read:
$$h_{i, j}(x, t)\,=\, (f_{i,j}(x),\, \gamma_{i, j}(t)),$$
where the $f_{i,j}$ do no depend on the $t$ variable and
are restrictions to their domains of definitions of translations in
$\RM^d$.

Note that each leaf of $\mathcal{L}$ is equipped with a
structure of oriented flat d-manifold inherited from the
Euclidean structure on $\RM^d$. Then the lamination satisfies the
further condition:
\smallskip

\noindent 2) Every leaf of $\mathcal{L}$, considered as an oriented flat 
$d$-manifold, is isometric to $\RM^d$.}
\end{defi}

\begin{defi}
\label{tilable}
{\rm A lamination $(M, \mathcal{L})$ is said {\it tilable} if it is flat and 
possesses a transversal $\Gamma$ which is a Cantor set.}
\end{defi}
\begin{defi}{\rm We call a domain $U$ of a chart $h:U\to V\times C$ of such a maximal atlas a {\it box} of the lamination. For any point $x$ in a box $B$  with coordinates $(p(x) , c(x))$ in the chart $h$, the slice $h^{-1} (V\times c(x))$ is called the horizontal of $x$ and the Cantor set $h^{-1}(\{p(x)\}\times C)$ is called the {\it vertical}  of $x$ in the box $B$.}\end{defi} 
Notice that these definitions make sense since the transition maps map  a vertical on  a vertical and  a horizontal on a horizontal. 

\noindent On a tilable lamination $(M, \mathcal{L})$ it is  possible to
define a parallel transport along each leaf of the lamination. We
denote by $Par (M, \mathcal{L})$ the set of parallel vector-fields
on $(M, \mathcal{L})$. This set is parameterized by $\RM^d$. This
gives a meaning to the notion of translation in the lamination by a
vector $u$ in $\RM^d$, defining this way a dynamical system $(M,
\omega_{\mathcal{L}})$. The lamination is {\it minimal} if such a
dynamical system is minimal.

Consider a perfect tiling $T\in T(X)$ constructed from a finite
collection of prototiles $X=\{p_1, \dots, p_k \}$. For each prototile $p_i$ in $X$, fix one interior point $y_i$. The set $Y =\{y_1, \dots , y_n\}$ determines a punctured version $X_Y$ of $X$. We denote by $\Omega_{T,Y}$  the canonical transversal associated to $Y$. From Proposition \ref{rectif}, we know that $\Omega_{T,Y}$ is a Cantor set.  It is plain to check that the dependence  on the set of marked points is continuous:
\begin{prop}
\label{Ycontinu}
Let us denote by $F(\Omega_T)$ the set of compact subsets of $\Omega_T$ endowed 
with the Hausdorff distance. Then the map $Y\to\Omega_{T,Y}$ from $\prod_j \, 
{\rm Int}(p_j)$ to $F(\Omega_T)$, is continuous.
\end{prop}

\noindent The topology of the Cantor set $\Omega_{T,Y}$ is generated by the
countable collection of {\it clopen} sets $U_{T,P}$ where $P$ is a
patch of a tiling in $\Omega_{T,Y}$ which contains the origin $0$
of $\RM^d$ and $U_{T,P}$ consists in these tilings in
$\Omega_{T,Y}$ whose restrictions to $P$ coincide with $P$.

\begin{teo}
\label{lam}
Let $T$ be a perfect tiling. Then there is a minimal tilable
lamination  $(\Omega_T,\mathcal{L})$  such that:

1) $(\Omega_T,\omega_{\mathcal{L}})$ and $(\Omega_T,\omega)$ are
conjugate dynamical systems.
\smallskip

2) For any arbitrary choice of the finite set  $Y=\cup_i
\{y_i\}$ where each point $y_i$ in the interior of each $p_i$'s, the set $\Omega_{T,Y}$ is a transversal
of the lamination.
\end{teo}
{\it Proof.} Fix one $\Omega_{T,Y}$ and choose at tiling $T'$ in $\Omega_{T,Y}$ and a patch $P'$ of $T'$ that contains in its interior,  a tile containing $0$. For each such pair	$(T',P')$ we consider the clopen set $C{(T',P')}$ in $\Omega_{T,Y}$ which consists in the tilings of $\Omega_{T,Y}$ that coincide with $T'$ on the interior of $P'$and we define 

$$U(T',P')\, =\, T"\, +\, v \,\,\,\forall \,\,\, (v,\, T")\in  int (P')\times C{(T',P')}.$$

We denote  $\phi(T', P')$ the map:

$$int (P')\times C{(T',P')}\to  U(T',P') $$

$$(v, \, T")\mapsto T"\, +\, v.$$

 When $T'$ runs over $\Omega_{T,Y}$ the $U(T', P')$'s form a cover of $\Omega_T$ from which we can extract a finite sub-cover, that we denote for sake of simplicity, $U_1, \dots , U_n,$ denoting also $P_1, \dots, P_n$ the corresponding patches , $C_1, \dots, C_n$ the corresponding clopen sets and $\phi_1,\dots , \phi_n$ the corresponding maps.  Each clopen set $C_i$ can be decomposed in a collection of pairwise disjoint clopen sets $C_{i, j}$, $j=1, \dots, k(i)$, with arbitrary small diameters.  Since the tilings in $\Omega_T$ are perfect, we can  choose these diameters small enough so that:

\begin{itemize}

\item For each $j=1, \dots , k(i)$ and each $i=1, \dots, n,$

the restriction of $\phi_i$ to $int(P_i)\times C_{i,j}$ is a homeomorphism onto its image that we denote by $U_{i,j}$;

\item For each pair $(i_1, i_2)$ in $\{1, \dots, n\}^2$ and for each $j_1$ in $\{1, \dots , k(i_1)\}$ there exists at most one $j_2$ in $\{1, \dots , k(i_2)\}$ such that $U_{i_1, j_1}\cap U_{i_2, j_2}\neq \emptyset$.  A tiling $\hat T$ in $\,U_{i_1, j_1}\cap U_{i_2, j_2}$ can be seen in two different ways: either as $T_1 +v_1 = \phi_{i_1}(v_1, T_1)$ or as $T_2 +v_2 = \phi_{i_2}(v_2, T_2)$ and it is important to notice that the vector $t_{1,2} = v_2-v_1$ is independent of the choice of $\hat T$ in $\,U_{i_1, j_1}\cap U_{i_1, j_1}$. The domain of definition of the transition map 

$\phi_{i_2}^{-1}\circ\phi_{i_1}$ is $int(P_{i_1})\cap (int (P_{i_2}) - t_{1, 2})\times C_{i_1, j_1}\cap (C_{i_2, j_2} - t_{1, 2})$ and on this domain, the transition map $\phi_{i_2}^{-1}\circ\phi_{i_1}(v_1, T_1)$ reads: 

$$\phi_{i_2}^{-1}\circ\phi_{i_1}(v_1, T_1)\, =\,(v_2, T_2)\,=\, (v_1 +t_{1,2}, T_1+t_{1,2}).$$

\end{itemize} 

 It follows that $\Omega_T$ can be endowed with a
structure of tilable lamination $\mathcal{L}$, where the leaf
containing a tiling $T'$ in $\Omega_T$ is  the $\RM^d$-orbit
$o(T')$ of $T'$. The set $\Omega_{T,Y}$ is a transversal of the lamination
by construction and it is plain to check that any other
$\Omega_{T,Y'}$ is also a transversal. A translation in the
lamination $(\Omega_T, \mathcal{L})$ defined using $Par (\Omega_T,
\mathcal{L})$ coincides by construction with a usual translation
acting on $\Omega_T$, so the last statement is immediate and
minimality is a direct consequence of the minimality of
$(\Omega_T,\omega)$.

\medskip
\begin{remark}
\label{fibreontorus}
{\rm Recently L. Sadun and R.F. Williams \cite{S-W}
have proved that the continuous hull of a perfect tiling of $\RM^d$ is 
homeomorphic to a bundle over the d-torus whose fiber is a Cantor set, a much nicer
object than a lamination. Unfortunately, this homeomorphism is not
a conjugacy: it does not commute with translations. However this shows that every 
tiling is orbit equivalent to a $\ZM^d$-action.}

\end{remark}


If at the first sight, tilings seem to be very rigid objects, their connections with laminations show that tiling spaces seen as  dynamical systems are easier to handle with. This 

 possibility of working with laminations will in turn lead to new results on tilings. For this purpose and by analogy with Markov partitions in Dynamical Systems, let us introduce the notion of {\it box decomposition } of a tilable lamination.


\begin{defi}{\rm A {\it well oriented $d$- rectangle} in $\RM^d$ is a $d$-rectangle  such that, for $0 {\leq} i \leq n$,  each $i$-face is parallel to an $i$-space generated by $i$ vectors of the canonical basis in $\RM^d$. A {\it block} in $\RM^d$ is a connected set  which is  a finite union of well oriented $d$-rectangles. A well oriented $d$-rectangle  is a {\it well oriented $d$-cube} when it is a $d$-cube.}\end{defi}


\begin{defi}\label{boite}{\rm  Consider a  tilable lamination $(M, \mathcal{L})$.

Consider a box $B\subset M$ which reads $int (P)\times  C$ in  a chart $h$ in $\L$  defined in a neighborhood of the closure $cl(B)$ where $C$ is a clopen subset of a Cantor set.

\begin{itemize}

\item when  $P$ is a well oriented $d$-rectangle, $B$ is called a box {\it  of rectangular type};

\item when  $P$ is a block  , $B$ is called a box {\it  of block type}.

\end{itemize}}

\end{defi}

Notice that these definitions are independent on the choice of the chart in $\L$.


\begin{defi}{\rm For any  box $B$ (of block or rectangular type) which reads $(int (P)\times C)$ in  a defining  chart, the set $h^{-1}(\partial P\times C)$ is called the {\it vertical boundary} of the box $B$. }\end{defi}


\begin{lemma} \label{rectangle} Consider a finite collection boxes of block type   $B_1, \dots, B_m$ of a tilable lamination. Then, there exists a finite collection of boxes of rectangular type in $B'_1\dots B'_p$  such that:

 \begin{itemize}

\item the $B'_l$'s are pairwise disjoint;

\item the closure of $\cup_{l=1}^{l=p}B'_l$ coincides with the closure of  $\cup_{j=1}^{j=m}B_j$;

\item if a $B'_l$ intersects a $B_j$ then it is contained in  this $B_j$.

\end{itemize}

\end{lemma}


\noindent{\it Proof}

\noindent Let us first prove this result with two boxes $B_1$ and $B_2$. For $i=1, 2$, let  $h_i$  be the chart defining $B_i$. We set $h_i(B_i) = int (P_i)\times C_i$. On its domain of definition, the transition map $h_{1, 2}$ reads:

$$h_{1, 2}(x, t)\,=\,h_1\circ h_2^{-1}\,=\, (f_{1,2}(x),\, \gamma_{1, 2}(t)),$$  

where $f_{1,2}$ does no depend on the $t$ variable and

is the  restriction  of a translation in

$\RM^d$ by a vector $t_{1, 2}$. Since, for $i= 1, 2$, there  exists a partition of the Cantor sets $C_i$ in clopen sets $C_{i, j}$, \, $j =1, \dots , k(i)$ with arbitrary small  diameters, we can choose such a partition so that $int (P_1)\times C_{1, j_1}\cap h_{1, 2}(int (P_2)\times C_{2, j_2})$ is the (possibly empty)  product space $int(P_1)\cap (int (P_2) - t_{1, 2})\times \tilde C_1$ where $\tilde C_1$ is a clopen subset of $C_{1, j_1}$. Similarly   $int (P_2)\times C_{2, j_2}\cap h_{2, 1}(int (P_1)\times C_{1, j_1})$ is the (possibly empty)  product space $int(P_2)\cap (int (P_1) + t_{1, 2})\times \tilde C_2$ where $\tilde C_2$ is a clopen subset of $C_{2, j_2}$.It follows that $h_1^{-1}( int (P_1)\times C_{1, j_1})\cup h_2^{-1}( int (P_2)\times C_{1, j_2})$ is the reunion of  five disjoint boxes of block type $ A_{1}\cup A_{2}\cup A_3\cup A_4\cup A_5,$

where

\begin{itemize}

\item $A_1\, =\, h_1^{-1}(int (P_1)\times C_{1,j_1}\setminus {\tilde C_1})$;

\item $A_2\, =\,  h_1^{-1}(int (P_1)\setminus (int (P_2) - t_{1, 2})\times  {\tilde C_1})$;

\item  $A_3\, =\, h_1^{-1}(int(P_1)\cap (int (P_2) - t_{1, 2})\times \tilde C_1)$;

\item $ A_4\,=\,  h_2 ^{-1}(int (P_2)\setminus (int (P_1) + t_{1, 2})\times  {\tilde C_2})$;

\item $A_5\, =\, h_2^{-1}(int (P_2)\times C_{2,j_2}\setminus {\tilde C_2})$.

\end{itemize}

Since $P_1,\,P_2,\,P_1\cap (P_2 - t_{1, 2}),\, P_1\setminus (P_2 - t_{1, 2})$ and $P_2\setminus (P_1 + t_{1, 2})$ are blocks, they are also  finite unions of well oriented $d$-rectangles $P'_1, \dots, P'_m$. This decomposition in $d$-rectangles induces a finite collection of disjoint rectangular boxes whose closure is  $h_1^{-1}( int (P_1)\times C_{1, j_1})\cup h_2^{-1}( int (P_2)\times C_{1, j_2})$ . Applying the same procedure for all pairs $(i_1, j_1)$ and $(i_2, j_2)$ we get the desired finite collection of rectangular boxes. The result for a collection of $n$ boxes in general is obtained by an easy induction.


\begin{defi} {\rm A tilable lamination $(M, \mathcal{L})$ admits a {\it box decomposition of rectangular type (resp. of block type)} if there exists a finite collection of boxes of rectangular type (resp. of block type) ${\cal B} =\{B_1, \dots , B_n\}$ such that:

\begin{itemize}

\item the $B_i$'s are pairwise disjoint;

\item the union of the closures of the $B_i$'s covers the whole lamination $M$.

\end{itemize}}

\end{defi}


\begin{prop} \label{boxon} Any tilable lamination $(M, \mathcal{L})$ admits a box decomposition of rectangular type.\end{prop}


\noindent {\it Proof.} Let $x$ be a point in $M$, $U_i$ an open set in $M$ containing $x$ and $h_i:U_i\to V_i\times T_i$ a chart. Choose a small well oriented $d$-rectangle $P_i$ in
$V_i$  and a small clopen set $C_i$ in $T_i$ such that $x\in U'_i = h_i^{-1}(int (P_i)\times C_i)$. Doing this
for all $x$ in $M$ we construct a cover of $M$ with the boxes of rectangular type $U'_i$. As $M$ is compact, we can extract form this cover a finite
one that we still denote $(U'_i)_{i=1, \dots, n}$. The proof is then a direct consequence of Lemma \ref{rectangle}.

\medskip

In terms of tilings, this last result implies that  the tilings induced on each leaf of
the lamination by the intersection of the leaf with the vertical boundaries
of the boxes $B_1, \dots , B_n$ defined in Proposition  \ref{boxon} are  made with at most $n$ $t$-types of well oriented $d$-rectangles. In the case when the lamination is minimal, they are perfect.
This is resumed in the next corollary which is  a converse to theorem \ref{lam} :


\begin{cor}
\label{til}
Consider a minimal tilable lamination $(M, \mathcal{L})$; then
there exists a  perfect tiling $T$ of $\RM^d$  made with a finite number of prototiles which are well oriented $d$-rectangles and such that
$\Omega_T = M$ (in the sense of theorem \ref{lam}).
\end{cor}


\medskip

\noindent
As an immediate consequence of  Theorem \ref{lam} and Corollary \ref{til} we have (see definition \ref{omegaequiv}):


\begin{cor} \label{cor-poly}
For any perfect tiling $T$ of $\RM^d$ there exists a  perfect
tiling $T'$   made with finite number of well oriented rectangular prototiles which is $\Omega$-equivalent to $T$.
\end{cor}


This should recommend us to work only with  perfect
tiling $T'$   made with a finite number  prototiles which are well oriented $d$-rectangles.  However we will perform construction on $\Omega$-equivalent tilings that will force us to leave this category. For this reason from now on we shall focus on perfect tilings made with a finite number of prototiles which are blocks and, unless otherwise stated, we will work with box decomposition of block type .


Notice that there is a lot of freedom in the construction  of the box decomposition of a tilable lamination done in Proposition \ref{boxon}.

\begin{defi}\label{sout}{\rm  Let ${\cal B}$  and ${\cal B}'$ be two box decompositions of block type of a same tilable lamination $(M, \mathcal{L})$. We say that  ${\cal B}'$ is {\it zoomed out} of ${\cal B}$ if:

\begin{itemize}

\item [(1)] for each point $x$ in a box $B'$ in ${\cal B}'$ and in  a box $B$ in ${\cal B}$, the vertical of $x$ in $B'$ is contained in the vertical of $x$ in $B$;

\item [(2)] the vertical boundaries of the boxes of ${\cal B'}$ are contained in the vertical boundaries of the boxes of ${\cal B}$;

\item [(3)] for each box $B'$ in ${\cal B}'$, there exists a box $B$ in ${\cal B}$ such that $B\cap B'\neq \emptyset $ and  the vertical boundary of $B$ does not intersect the vertical boundary of $B'$. 

 \end{itemize}}

\end{defi}

\begin{prop}\label{zoomout}Consider a minimal tilable lamination $(M, \mathcal{L})$. Then, for any box decomposition of block type ${\cal B}$, there exists another box decomposition of block type zoomed out of ${\cal B}$.\end{prop}

\noindent{\it Proof.}

\noindent{\it Step 0.} Let ${\cal B} =\{B_1, \dots , B_n\}$ be a box decomposition of  block type of $M$ and for $i= 1, \dots, n$, let $h_i$ be the chart mapping $B_i$ on $P_i\times C_i$. As already observed in the proof of lemma \ref{rectangle}, up to a splitting of  the vertical of the point in each box, we can  assume that the intersection $B_i\cap B_j$ for $i, j$ in $\{1, \dots , n\}$   (read in a chart) is the product (possibly empty) of a block by a clopen set.

\noindent{\it Step 1.}  For each $i=1, \dots, n$, we  denote by $J(i)$ a set of $l$'s in $\{1, \dots , m\}$ such that $cl(B_l)\cap cl(B_i)\neq \emptyset,$ and define:
$$W_i\,=\, \cup _{l\in J(i)} cl(B_l).$$ 
Let $x$ be in $B_i$ and consider in $W_i$ the connected component of leaf of the lamination containing $x$, say $L_x$. For any subset $K$ of $J(i)$ we denote by $B_{i, K}$ the box made of all these points $x$ in $B_i$ such that $L_x$ intersects $B_l$ if and only if $l$ is in $K$. We get a partition:
$$B_i\, =\, \cup_{K\subset J(i)}B_{i, K}.$$

For each $K$ in $J(i)$, we extend the box $B_{i, K}$ as follows:

$$B^{(1)}_{i, K}\, =\, \cup_{x\in B_{i, K}}int (L_x).$$

If the vertical sizes of the boxes of the box decomposition $\cal B$ are chosen small enough and since the leaves of the lamination are copies of $\RM^d$,  the open sets $B^{(1)}_{i, K}$ are boxes    and for each $i=1, \dots, n$, two boxes $B^{(1)}_{i, K}$ and $B^{(1)}_{i, K'}$ are disjoint when $K$ and $K'$ are disjoint  subsets of $J(i)$.
Consider the collection of boxes $B^{(1)}_{1, K}$ when $K$
runs over the subsets of $J(1)$ and,  to get simpler notations,  call this sequence $B^{(1)}_1, \dots, B^{(1)}_m$.

\noindent For $l=2, \dots, n$, consider  the boxes $B^{(1)}_{m+l-1} = B_l\setminus (\cup_{k=1}^{m+1} cl(B^{(1)}_k))$. It is clear that the sequence of boxes ${\cal B}^{(1)} =\{B^{(1)}_1, \dots , B^{(1)}_{m+n-1}\}$ is a box decomposition of the lamination $M$.  This box decomposition satisfies points $(1)$ and $(2)$ of Definition \ref{zoomout}. It is not (a priori) zoomed out of ${\cal B}$ since point $(3)$ is satisfied only by the boxes $B^{(1)}_1, \dots, B^{(1)}_m$ which cover $B_1$.


\noindent{\it Step 2.} We Consider now the box decomposition ${\cal B}^{(1)} =\{B^{(1)}_1, \dots , B^{(1)}_{m+n-1}\}$   and perform to the box $B^{(1)}_{m+1}$ the same procedure as the one we did in Step 1 for $B_1$. We get this way a third box decomposition,  ${\cal B}^{(2)} =\{B^{(2)}_1, \dots , B^{(2)}_{n_2}\}$  of the lamination $M$.  This box decomposition satisfies points $(1)$ and $(2)$ of Definition \ref{zoomout} and point $(3)$ is satisfied only by the boxes $B^{(1)}_1, \dots, B^{(1)}_m$ which cover $B_1$ and $B_2$.


\noindent{Step 3.} We can iterate this procedure to get after $n$ steps, a  box decomposition ${\cal B}^{n} = \{B^{(n)}_1, \dots , B^{(n)}_{m_n}\}$ zoomed out of ${\cal B}$.

\medskip
Consider two box decompositions ${\cal B}= \{B_1, \dots, B_n\}$ and ${\cal B}'=\{B'_1, \dots, B'_m\}$ of a same  tilable lamination $(M, \mathcal{L})$ such that ${\cal B}'$ is zoomed out of ${\cal B}$. For each $j$ in $\{1, \dots , m\}$, we denote by $I(j)$ the sets of $i$'s in $\{1, \dots , m\}$ such that $cl(B_i)\cap cl (B_j)\neq \emptyset$ and $B_i\cap B'_j = \emptyset$ and we define:

$$Z_j\, = \cup_{i\in I(j)}cl(B_i).$$

\begin{defi}{\rm  We say that ${\cal B}'$ {\it forces its border}\footnote{Actually the concept of "forcing its border" was first introduced in the context of tilings (see \cite{KP00}). We will see in the sequel how our definition fits with the standard one.}
 if for  each $j$ in $\{1, \dots , m\}$ and for each $x$ in $B'_j$, the connected component $L_x$ in $Z_j$ of  the leaf of the lamination containing $x$ intersects $B_i\setminus B'_j$ if  and only if  $i$ is in $I(j)$.}\end{defi}

Actually Proposition \ref{zoomout} can be improved in the following theorem:

\begin{teo}{\label{zoomoutforce}}Consider a minimal tilable lamination $(M, \mathcal{L})$. Then, for any box decomposition of block type ${\cal B}$, there exists another box decomposition of block type zoomed out of ${\cal B}$ that forces its border.\end{teo}

\noindent{\it Proof.}
Fix $j$ in $\{1, \dots , m\}$, and for each subset $H$ of $I(j)$ we  denote by $B'_{j, H}$ the box made of all these points $x$ in $B'_j$ such that $L_x$ intersects $B_l\setminus B'_j$ if and only if $l$ is in $H$. We get a partition:

$$B'_j\, =\, \cup_{H\subset I(j)}B'_{j, H}.$$

It is clear that the box decomposition made of all the $B'_{j, H}$'s, $H\subset I(j)$ and $j$ in $\{1, \dots , m\}$ is a box decomposition which is zoomed out of ${\cal B}$ and forces its border. 


\begin{cor}\label{corect}
Consider a minimal tilable lamination $(M, \mathcal{L})$. Then, for any box decomposition of block type ${\cal B}$, there exists a sequence of  box decompositions of block type ${\cal B}^{(n)}_{n\geq 0}$ such that:

\begin{itemize}

\item ${\cal B}^{(0)} = {\cal B}$;

\item for each $n\geq 0$, ${\cal B}^{(n+1)}$ is zoomed out of ${\cal B}^{(n)}$ and forces its border.\end{itemize}

\end{cor}

In terms of tilings, this last result can be interpreted as follows. Consider a leaf $L$ of the lamination, fix $n\geq 0$, and consider the  tiling $T^{(n)}$ on this leaf induced  by the intersection of the leaf with the vertical boundaries of the boxes of the box decomposition ${\cal B}^{(n)}$ defined in Theorem  \ref{zoomoutforce}. 

\begin{defi}\label{nested}{\rm A sequence $T^{(n)}_{n\geq 0}$ of tilings obtained this way is called {\it a nested sequence of tilings}.}
\end{defi}

It satisfies in particular the following properties:

\noindent For all $n\geq 0$:

\begin{itemize}

\item the $t$-tiles types of $T^{(n)}$ coincide with a finite number of prototiles which are blocks;

\item each  tile $t^{(n+1)}$ of $T^{(n+1)}$ is a connected patch of $P^{(n)}$ of  $T^{(n)}$ and contains at least one tile of $T^{(n)}$ its interior;

\item furthermore, the patch $Q^{(n)}$ of $T^{(n)}$ made with all the tiles of $T^{(n)}$ that intersects $P^{(n)}$  is uniquely determined up to translation\footnote{this is a generalization of the classical definition of "forcing its border" for substitution  tilings (see \cite{KP00}).};

\item the tiling  $T^{(n)}$ is $\Omega$-conjugate to $T^{(0)}$.\end{itemize}

Consider a minimal tilable lamination $(M, \mathcal{L})$ and  a box decomposition ${\cal B}$. We define an equivalence relation on $M$ as follows: two points in $M$ are equivalent if the are in a same box of  ${\cal B}$ and on the same vertical in this box. The  quotient space is a branched $d$-manifold that inherits from its very construction some extra structures. Section  \ref{sec-BOFS} is devoted to an axiomatic approach of these objects. If we perform this quotient operation for a sequence of box decompositions zoomed out one of the other as in Corollary \ref{corect}, we get a sequence of such branched manifolds whose study is reported in Section \ref{sec-tilingEFS}.

\section{BOF-d-manifolds}
\label{sec-BOFS}

\noindent In this section we will establish the basic notions concerning the
category of {\it branched oriented flat} (briefly: BOF) $d$-manifolds.
There exist several possible variations on these notions. We will try
to adopt the simplest version suitable to  applications to tilings. 

\medskip

\subsection {Local Models}
\label{ssec-local}

For $r> 0$, and $x= (x_1, \dots , x_d)$ in $\RM^d$,  we denote by ${\cal C}(r, x)$ the $d$-dimensional open cube in $\RM^d$ centered at $x$ with size $r$:

$${\cal C}(r, x) \,=\{(y_1, \dots , y_d)\,\, \vert \,\, sup_{i=1, \dots , d}\vert y_i-x_i\vert < r\}.$$

For $i = 1, \dots , d$, let $\epsilon_i =\pm 1$ and consider the $d$-orthant 
$ O_{\epsilon_1, \dots , \epsilon_d}(r,x)$ which consists in the closure of the points $y$ in $ {\cal C}(r,x)$ such that the sign of $y_i-x_i$ is $\epsilon _i$ for $i=1, \dots , d$. 

\begin{defi}
\label{d-cube}
{\rm 

\begin{itemize}

\item A {\it type 1 $d$-cube} (with size $r>0$, and centered at $x$) is the  $d$-cube $({\cal C}(r, x))$.
\item For $1\leq p\leq 2^d$, a {\it  type $p\,$ $d$-cube}  (with size $r>0$, and centered at $x$) is given by the data
$({\cal C}(r,x), ({\cal P}_1,  n_1), \dots ({\cal P}_p, n_p))$ where:

\begin{itemize}

\item [(i)]   ${\cal P}_1\dots , {\cal P}_p$ are $p$ collections of $d$-orthants ${\cal P}_1\dots , {\cal P}_p$ such that \-

$ int ({\cal P}_i)\cap int ({\cal P}_j) = \emptyset$ for $i\neq j$ and $\cup_{i=1}^{i=p}{\cal P}_i = {\cal C}(r)$;

\item [(ii)] for $i=1, \dots p$, $n_i$ is a  positive integer associated to ${\cal P}_i$.

\end{itemize}

Consider the set of points ${\cal X} =\{(y, n)\}$ where $y$ is in ${\cal C}(r,x)$ and $n=n_i$ whenever $y$ is in ${\cal P}_i$. 
The corresponding  type $p\,$ $d$-cube (of size $r>0$ centered at $x$) is the quotient space  ${\cal X}/\sim$
where $(y,  n)\sim (y', n')$ if and only if:

\begin{itemize}
\item $ (y,  n)\,=\, (y', n')$;
\item or $y\, =\, y'$ and belongs to at least two distinct collections of $d$-orthants  ${\cal P}_i$ and ${\cal P}_j$.
\end{itemize}

\item A {BOF-$d$-cube} is a type $p\,$ $d$-cube  with some  size $r>0$ and centered at some
 $x$ in $RM^d$, for some $p$ in $1, \dots 2^d$.\end{itemize}
}
\end{defi}

Consider a BOF-$d$-cube $\mathcal{D}$.  

\noindent $\bullet$  The tangent plane $T_{\tilde x}\mathcal{D}$ to $\mathcal{D}$ at every
point $\tilde x$ in $\mathcal{D}$ is well defined; thus the
BOF-$d$-cube $\mathcal{D}$  has a natural {\it branched
$C^1$-structure}. In fact the {\it tangent bundle}
$T(\mathcal{D})=\cup_{\tilde x}T_{\tilde x}\mathcal{D} $ is trivial: one
trivialization is associated with each positive basis of $\mr^d$. We
call {\it canonical trivialization} the trivialization associated with the
standard basis of $\mr^d$. In this way we have also specified one
{\it orientation} of $\mathcal{D}$.

\smallskip

\noindent $\bullet$ The standard Euclidean metric on $\mr^d$ ($ds^2=
dx_1^2+\dots + dx_d^2$) induces, in a natural way, a {\it branched flat
metric} on $\mathcal{D}$, whence a distance, a notion of
``branched geodesic arc'' and so on.

\smallskip

 \noindent $\bullet$ Each BOF-$d$-cube (with radius $r>0$ and centered at $x$) is the union of a finite copies of  a  well oriented $d$-cube ${\cal C}(r, x)$ glued along collections of $d$-orthants,  such that the orientation and both the branched $C^1$-structure and the branched flat metric restrict to the usual ones. Any such an embedded well oriented $d$-cube ${\cal C}(r, x)$ is called a {\it smooth sheet} of the  BOF-$d$-cube $\mathcal{D}$.

\noindent A type $p\,$ $d$-cube $\mathcal D$ is stratified as follows:
for $1\leq l\leq d$, $ {\cal V}_l $  is the set of points $(x, n)$ in 
$\mathcal D$ such that $x$ belongs to exactly $d-l+1$ distinct ${\cal P}_i$'s and  $ {\cal V}_0 $  is the set of points $(x, n)$ in 
$\mathcal D$ such that $x$ belongs to strictly more than $d$ distinct ${\cal P}_i$'s. We have:

\begin{itemize}

\item $\mathcal D\,= \, cl({\cal V}_d) \, =\,\cup_{l=0}^{l=d} {\cal V}_l$;

\item  for each $0< \leq d$, $\,cl({\cal V}_l)\, =\, {\cal V}_l\cup {\cal V}_{l-1}$ and ${\cal V}_l\cap {\cal V}_{l-1}=\emptyset$;

\item $cl({\cal V}_0) = {\cal V}_0$ is reduced to the center of $\mathcal D$.

\end{itemize}

For $1\leq l\leq d,$ each set $ cl({\cal V}_l)$ is a connected set  which is a collection of $l$-faces of well oriented $d$-cubes containing the center of $\mathcal D$. Furthermore each point in ${\cal V}_l$ is the center of a type $d-l+1\,$ $d$-cube in $\mathcal D$.

\begin{defi}{\rm For $0\leq l\leq d$ the set$ {\cal V}_l$ is called the {\it $l$-stratum} of the  BOF-$d$-disk $\mathcal D$.}
\end{defi}.

Let us define now  adapted morphisms between the local models.

\begin{defi}
\label{localsub}
{\rm A continuous map $f:\mathcal{D} \to \mathcal{D}'$ between two
BOF-$d$-cubes of the same size $r>0$ is  a {\it local BOF-submersion}
(onto its image) if:

1) $f$ is $C^1$ with respect to the branched $C^1$-structures.

\smallskip
2) For every open smooth sheet $D$ of $\mathcal{D}$, $f(D)$ is a
smooth sheet of $\mathcal{D}'$; the restriction of $f$ to $D$
coincides with the restriction of a translation. In fact, there is
one translation which works for all smooth sheets of $\mathcal{D}$. 
 If one fixes
one basis of $\mr^d$ and takes the corresponding trivializations
of the two tangent bundles, then the differential $df_x$ of $f$ at
any point $x$ of $\mathcal{D}$, induces the identity on $\mr^d$
(which has been identified with both $T_x\mathcal{D}$ and
$T_{f(x)}\mathcal{D}'$). }
\end{defi}

\begin{defi}
\label{isometry}
{\rm A local BOF-submersion $f$ is  a {\it local 
BOF-isometry} iff it is bijective. In such a case $f(\mathcal{D})$ is a BOF-$d$-cube and $f^{-1}$ is a local BOF-isometry and the stratification associated with $\mathcal D$ is mapped on the stratification associated with $f(\mathcal{D})$ (stratum by stratum).}
\end{defi}

\subsection {\bf BOF-$d$-manifolds}
\label{ssec-bofsurf}

We are now in a position to extend our construction to global models.
\begin{defi}
\label{bofs}
{\rm A {\it BOF-$d$-manifold} $B$ is a compact,
connected  metrizable topological space endowed with a (maximal) atlas $\{U_j,\phi_j\}$ such that:

1) every $\phi_j: U_j\to W_j$ is a homeomorphism onto an open set
$W_j$ of some open BOF-$d$-cube.

2) For any BOF-$d$-cube $\mathcal{D}$ embedded into $\phi_j(W_i\cap
W_j)$, the restriction of $\phi_{ij}=\phi_i\circ \phi_j^{-1}$ to
$\mathcal{D}$ is a local BOF-isometry onto a BOF-$d$-cube embedded
into $\phi_i(W_i\cap W_j)$.}
\end{defi}

Hence for every point $x$ of a BOF-$d$-manifold $B$, there exists a neighborhood $U$ of $x$ and a chart $\phi: U\to \phi (U)$ such that $\phi(U)$ is a BOF-$d$-cube and $\phi(x) =0$. The type of the BOF-$d$-cube $\phi(U)$ is uniquely determined by $x$. Such a neighborhood $U$ is called a{\it \, normal
neighborhood} of $x$.

\begin{defi}\label{injrad}{\rm
The {\it injectivity radius of $x\in B$}, denoted by inj$_B(x)$,
is the sup of the size of any normal neighborhood of $x$ in $B$. The {\it injectivity radius of $B$}, denoted $inj(B)$ is defined by:
$$ {\rm inj}(B) = \inf_{x\in B} {\rm inj}_B(x)\ .$$}
\end{defi}

\begin{defi}\label{natstrat} {\rm For $0\leq l\leq d$, the {\it  $l$-face} of a BOF-$d$-manifold $B$ is the set of points $x$ in $B$ such that there exits  a chart $(U, \phi)$ in the atlas such that $x$ is in $U$ and $\phi(x)$  belongs to the $l$-stratum of the BOF-$d$-cube $\phi (U)$. Notice that this property is independent on the choice of the chart $(U, \phi)$ as long as $x$ is in $U$. A {\it $l$-region} is a connected component of the $l$-face of $B$.  The finite partition of $B$ in $l$-region, for $0\leq l\leq d$, is called the {\it natural stratification} of $B$. The union of all the $l$-regions for $0\leq l\leq d-1$,  forms the {\it singular locus}
$Sing(B)$ of $B$.}
\end{defi}

Note that, in particular, any oriented flat $d$-torus $B$ is a BOF-$d$-manifold with one $d$-region $R=B$. In general, any $d$-region of a BOF-$d$-manifold $B$ is a
connected (in general non compact) $d$-manifold naturally endowed with
a $(X,G)=(\mr^d,\mr^d)$-structure, where the group $G=\mr^d$ acts
on $X=\mr^d$ by translation. This means that any region admits a
$d$-manifold atlas such that all the transition maps on connected
domains are restriction of translations. One can see for
instance  \cite{B-P} for the basic notions about the
$(X,G)$-structures set up.

\begin{remark}
\label{bofsrem} 
{\rm All the objects that we have associated with any open BOF-$d$-cube (the tangent 
bundle, its trivializations, the orientation, the branched $C^1$-structure, the 
branched flat metric and so on) globalize to any BOF-$d$-manifold. In particular, 
one has a natural notion of parallel transport on $B$ with respect to the flat 
metric, whence the notion of {\it parallel vector field} on $B$. Let us denote by 
$Par(B)$ the set of parallel vector fields on $B$. Fix the canonical 
trivialization of the tangent bundle $T(B)$, then we have a natural isomorphism
$$\rho_B: \mr^d \to Par(B)$$

\noindent defined by identifying $\mr^d$ with any tangent space $T_{\tilde x}B$, and by 
associating to every  vector $v\in \mr^d$ the
parallel vector field $\rho_B(v)$ obtained by parallel transport
of $v$, starting from $\tilde x$. This does not depend on the choice of
the base point $\tilde x$.}
\end{remark}

\begin{defi}\label{cellbofs}{\rm A  BOF-$d$-manifold $B$  has {\it rectangular faces}  if, for   each  $d$-region $R$ of $B$,  there exists a $C^1$ injective map $f: R\to \RM^d$ such that $f(R)$ is the interior of a well oriented $d$-rectangle  in $\RM^d$ and the differential of $f$, read in the charts of $B$, satisfies $df_x =Id$ at each point $x$ in $R$.}
\end{defi}

In particular, the $d$-regions of a BOF-$d$-manifolds do not carry any topology.

\noindent {\bf Convention.} From now on, we will only be concerned with  BOF-$d$-manifolds with rectangular  faces (unless otherwise stated).

\subsection {BOF-submersion}
\label{ssec-bofsubmer}
\begin{defi}\label{defbofsub}{\rm A continuous map $f: B\to B'$
between BOF-$d$-manifolds is a  {\it BOF-submersion} if:
\smallskip

1) $f$ is $C^1$ and surjective.

\smallskip

2) For every $x\in B$ and for every normal neighborhood $\mathcal{D}'$ of $f(x)$ in $B'$, there exists a normal neighborhood
$\mathcal{D}$ of $x$ in $B$ such that $f(\mathcal{D})\subset \mathcal{D}'$ and, read in the corresponding charts, 
$f|: \mathcal{D}\to f(\mathcal{D})$ is a local BOF-submersion.

\smallskip

3) For each region $R$ of $B$, there exists a region $R'$ of $B'$ such that $f$ is a diffeomorphism from $R'$ to $R$. In particular, the singular set of $B$ is mapped into the singular set of $B'$:
$$ Sing(B) \subset f^{-1}(Sing(B')).$$
}
\end{defi}

\begin{remark}\label{rayon} Notice that the pre-image of a normal neighborhood  $U'$ in $B'$ (with radius $r>0$) is a finite union of disjoint normal neighborhoods $U_1, \dots, U_n$ in $B$ (with radius $r>0$) and that $f: U_i\to U'$ is a $ C^1$ bijection whose differential is the identity when read in the charts. This imply in particular that:

$$ inj(B')\,\, \leq \,\, inj (B).$$            
\end{remark}

\subsection {Cycles and Positive Weight Systems}
\label{ssec-cyclepws}

\begin{defi}\label{pws}{\rm A {\it (strictly) positive weight system} $\,w$
on a BOF-$d$-manifold $B$ is a function which assigns to each $d$-region $R$ of
$B$ a real number ($w(R)>0$) $w(R)\geq 0$ in such a way that the
``switching rules'' (or Kirchoff-like laws) are satisfied; this
means that along every $d-1$-region $e$ of $B$ the sum of the weights on
the germs of $d$-regions along $e$ on one side equal the sum of the
weights of the germs of region on the other side. Let us denote by
($W^*(B)$) $W(B)$ the set of these weight systems. The BOF-$d$-manifold $B$ equipped
with $w\in W^*(B)$ is called a {\it measured} BOF-$d$-manifold. The {\it total
mass} $m(w)$ of $w\in W(B)$ is just the sum of all the weights
of $w$. So we have a partition $\{W^m(B)\}_{m\in \mr^+}$ of $W(B)$
by the different total masses. }
\end{defi}

\begin{remark}\label{rempws}{\rm
 Note that the definition of positive a weight system does not
involve any region orientation, and in fact it makes sense even
for non-orientable branched manifolds (see for instance \cite{B-P2}
for more details on this notion).}
\end{remark}

Let us now  exploit the fact that the $d$-region orientations are
part of the BOF-$d$-manifold structure.
Let us  first recall few elementary facts
about the (cellular) (co)-homology of $B$. 
Consider the natural stratification  of $B$, ${\cal V}_0, \dots {\cal V}_l$, where, for $i=0, \dots d$, ${\cal V}_i$ is decomposed into a finite number $a_i$ of $i$-region that we orient and order in an arbitrary way but for the orientation of the $d$-regions which is the natural orientation induced by the BOF $d$-manifold structure.
We set  $\ma = \mz,\, \mr \,$ and denote by $C_i(B,\ma)$  the free
$\ma$-module which has as (ordered) basis the set of ordered and
oriented $i$-regions   in ${\cal V}_i$. By convention, for any oriented $i$-region $e$, $-e = -1e$, and it consists of the same cell with the
opposite orientation. We call $C_i(B,\ma)$ the module of the
 $i$-{\it chains} of $B$. Clearly, as we dispose by
definition of a distinguished basis of $C_i(B,\ma)$, it can be
identified with $\ma^{a_i}$.

\begin{remark}\label{order}{\rm With the exception of the $d$-regions orientation,
there are no canonical choices for the above orders and orientations.
On the other hand, it is clear that two different choices reflect in a very
simple linear automorphism of $\ma^{a_i}$ which is completely under control.
Moreover, these choices will be essentially immaterial in our discussion.}
\end{remark}

We define the  linear {\it boundary operator}
$$\partial_{i+1} : C_{i+1}(B,\ma)\to  C_i(B,\ma)\ $$
which assign to any $i+1$-region, the sum of the $i$-regions that are in its closure pondered with a positive sign (resp. negative) if the induced orientation fits (resp. does not fit) with the orientation chosen for these $i$-regions. It is clear that $\partial_i\circ\partial_{i+1}=0$. 

\noindent The  space $Z_i(B,\ma)= {\rm Ker}\ \partial_i$ is called the space
of  $i$-{cycles} of $B$ and the space  $B_i(B,\ma)=
\partial_{i+1}(C_{i+1}(B,\ma))$ is called the space of
$i$-{boundaries} of $B$. In fact, $B_i(B,\ma)\subset Z_i(B,\ma)$
and $H_i(B,\ma)= Z_i(B,\ma)/ B_i(B,\ma)$ is the
$i^{th}$ homology group of $B$. Note that $H_d(B,\ma)=
Z_d(B,\ma)$. 

A standard result of algebraic topology insures that (up to $\ma$-module isomorphism) $H_i(B,\ma)$ is a
topological invariant of $B$ that coincides with the $i^{th}$ singular
homology of $B$ (see for example \cite{Spa}). It is important to observe  that a $d$-chain $z$ of $B$
is  a $d$-cycle (i.e. $\partial_d(z)=0$) if and only if the
coefficients of $z$ formally satisfy the above ``switching rules''
(or Kirchoff-like laws), extended to possibly non positive
``weights''. Note that here we are using the distinguished choice
of $d$-region orientations.
Let us denote by ($Z^{>0}_d(B,\ma)$) $Z^{\geq 0}_d(B,\ma)$ the
{\it (strict) positive cone} of $Z_d(B,\ma)$, which is formed by
the $d$-cycles with (strictly) positive coefficients. One clearly
has

\begin{lem}\label{W=Z} $W(B)=Z^{\geq 0}_d(B,\mr)$, $W^*(B)=Z^{>0}_d(B,\mr)$.
\end{lem}

 The dual module $C^i(B,\ma)$ of $C_i(B,\ma)$ is called the module of
$i$-co-chain of $B$. It can be identified with $\ma^{a_i}$ using
the dual basis of the above distinguished basis of $C_i(B,\ma)$.
There are linear {\it co-boundary} operators $$\delta^i:
C^i(B,\ma)\to C^{i+1}(B,\ma)$$ defined by the relation
$\delta^i(c)(z)=c(\partial_i(z))$. So we have the module of
co-cycles $Z^i(B,\ma) = {\rm Ker}\delta^i$, the module of
co-boundaries $B^i(B,\ma)=\delta^{i-1}(C^{i-1}(B,\ma))$, and
finally (up to module isomorphism) 
$$H^i(B,\ma)\,=\, Z^i(B,\ma)/B^i(B,\ma)\ .$$
Note that $Z^d(B,\ma)=C^d(B,\ma)$
which is the free $\ma$-module spanned by the {\it
characteristic functions} of the $d$-regions of $B$.  

If $\alpha \in  H^d(B,\ma)$ is
represented by a   $d$-co-cycle $c= \sum_i s_ie_i^*$ (remind
that we are using the dual basis) , and $w=\sum_i w_ie_i$ is a $d$-chain, then
$$\alpha(w)\, =\, <w\vert c>\,   = \,\sum_i s_i\omega_i \ .$$

\begin{remark}\label{misure-cochain}{\rm An immediate but relevant
 consequence is that one can
{\it contract any $d$-co-homology class against a ``positive
measure''} on $B$ (a fact which is not possible, for example, if
the branched manifold is not orientable). For any $w\in W(B)$, 
it is not necessary to know that
$w$ is cycle in order to define $<w|c>$, by using the last formula
above. It should make sense even for a non-orientable branched
manifold. What is remarkable in our situation is that
$<w|c>=<w|[c]>$, where $[c]$ is the co-homology class represented
by $c$. It follows that the natural pairing $< \vert >$ induces a pairing
$$<|>: W(B)\times H^d(B,\ma)\to \mr \ \ \ \ \ <w|\alpha>= \alpha (w).$$
}
\end{remark}

The following proposition is immediate.

\begin{prop}\label{subhom} If $f: B\to B'$ is a  BOF-submersion
then:

\smallskip

1) there exists a natural linear map (well defined up to the mild
ambiguity indicated in remark \ref{order}) $f_*: Z_d(B,\ma)\to
Z_d(B',\ma)$ such that $f_*(W^m(B))\subset W^m(B')$,
$f_*(W^{*,m}(B))\subset W^{*,m}(B')$, for every $m\in \mr^+$ .

\smallskip

2) There exists a natural linear map $f^*: C^d(B',\ma)\to
C^d(B,\ma)$ such that $[\alpha](f_*(w))= [f^*(\alpha)](w)$, for
every $\alpha \in  C^d(B',\ma)$, every $w\in Z_d(B,\ma)$, where
$[\alpha]\in H^d(B,\ma)$ is the co-homology class represented
by the co-chain $\alpha$. In fact (with a slight abuse of
notation), $f^*([\alpha])=[f^*(\alpha)]$ well-defines a linear map
$f^*: H^d(B',\ma)\to H^d(d,\ma)$.
\end{prop}

\begin{lem}\label{integral} Let $B$ be a BOF-$d$-manifold and
$\mu=(\mu_i,\dots,\mu_k)\in W(B)$. Let $c\in Z^d(B,\mr)$ which
represents a class of $H^d(B,\mz)$. Then there exists a sequence
of integers $(m_1,\dots,m_k)$ such that
$$ <\mu|c> = \sum_i m_i\mu_i \ .$$

Moreover, all such linear combinations with integer coefficients
arise in this way.
\end{lem}

{\it Proof.} As $\mu$ is a $d$-cycle the value of $<\mu|c>$ does not
change if one add to $c$ any co-boundary. As $c$ represents a
class in $H^d(B,\mz)$, $c$ differs from a suitable $c'\in
Z^d(B,\mz)$ by a co-boundary. The first statement of the lemma follows.
The last one is a consequence of the fact that every $d$-co-chain in
$C^d(B,\mz)$ is actually a $d$-co-cycle.

\section{Tilings \& Expanding Flattening Sequences}
\label{sec-tilingEFS}

\subsection{Expanding Flattening Sequences}
\label{ssec-EFS}

\begin{defi}
{\rm A  BOF-submersion $f: B\to B'$ satisfies
the {\it flattening condition}, if for every $x\in B$ and for every normal neighborhood $\mathcal{D}'$ of $f(x)$, there exists a  small enough normal
neighborhood $\mathcal{D}$ of $x$ in $B$ such that $f(\mathcal{D})\subset \mathcal{D}'$ and, read in corresponding charts, $f|: \mathcal{D}\to f(\mathcal{D})$ is a local BOF-submersion that maps $\mathcal {D}$ on one single sheet of   of $f(\mathcal{D})$.}
\end{defi}

\begin{defi}\label{defi-efs}
{\rm An {\it expanding flattening sequence} EFS is a
sequence $\mathcal{F}=\{ f_i:  B_{i+1}\to B_i\ \}_{i\in \mn}$ of submersions such that:
\begin{itemize}

\item [(i)] the sequence of injectivity radius of the $B_i$'s is a strictly increasing sequence that goes to $+\infty$ with $i$;

\item [(ii)] for each ${i\in \mn}$ the map $f_i$ satisfies the flattening condition.

\end{itemize}
}
\end{defi}

Once an EFS $\mathcal{F}$ is chosen, comes,  naturally associated with this sequence,  the further ``inverse'' or ``direct'' sequences
($\ma = \mz,\ \mr$):

i) $$Z_d(\mathcal{F},\ma) = \{ (f_i)_*:  Z_d(B_{i+1},\ma)\to
Z_d(B_i,\ma)\ \}_{i\in \mn}$$ with the restricted sequences
$$W(\mathcal{F},\ma) = \{ (f_i)_*: W(B_{i+1},\ma)\to W(B_i,\ma)\
\}_{i\in \mn}$$ and
$$W^*(\mathcal{F},\ma) = \{ (f_i)_*: W^*(B_{i+1},\ma)\to
W^*(B_i,\ma)\ \}_{i\in \mn}\ .$$

ii) $$C^d(\mathcal{F},\ma) = \{ (f_i)^*:  C^d(B_i,\ma)\to
C^d(B_{i+1},\ma)\ \}_{i\in \mn}$$ which induces
$$H^d(\mathcal{F},\ma) = \{ (f_i)^*:  H^d(B_i,\ma)\to
H^d(B_{i+1},\ma)\ \}_{i\in \mn}\ .$$

iii) $${\rm Par}\ (\mathcal{F}) = \{ df_i: {\rm Par}\ (B_{i+1})\to {\rm
Par}\ (B_i)\ \}_{i\in \mn}\ .$$

Associated with these "inverse " (resp. "direct" sequences, come their {\it inverse limits} (resp. the direct limits) which will be relevant in the next sections and whose definitions are briefly recalled here.

 Given an ``inverse'' sequence of maps $\mathcal{X}=\{X_i\stackrel{\tau_i}{\leftarrow} X_{i+1}\}_{i\in \mn}$, let us recall that the elements of the {\it projective limit set} (or {\it inverse} limit)
$$\lim_\leftarrow \mathcal{X}$$
 consists of the elements  $(x_0, x_1, \dots, x_n, \dots)$ in the product $\prod_{i\geq 0} X_i$ such
that $\tau_i(x_{i+1}) = x_i$ for all $i\geq 0$.

\noindent Note that for every $j\geq 0$ there exists a natural map
$p_j: \lim_\leftarrow \mathcal{X} \to X_j$. If the $X_j$'s are
topological spaces and the maps continuous maps, the set
$\lim_\leftarrow \mathcal{X}$ is a topological space with the
finest topology such that all the $p_j$'s are continuous. Hence it
is a subspace of the topological product $\prod_{i\geq 0} X_i$. Thus, if all the $X_i$'s are compact, $\lim_\leftarrow \mathcal{X}$ is compact.

 Similarly, given any ``direct'' sequence of maps
$\mathcal{Y}=\{ Y_i\stackrel{\tau_i}{\rightarrow} Y_{i+1} \}$, the
elements of the {\it direct limit set} 
$$\lim_\rightarrow \mathcal{Y}$$ 
are the sequence of the form $(x_j, x_{j+1}, \dots, x_n, \dots)$, for some $j\geq 0$ such that $x_h\in Y_h$ and
$\tau_h(x_{h}) = x_{h+1}$. Here again note that for every $j\geq
0$ there exists a natural map $i_j: Y_j \to \lim_\rightarrow
\mathcal{Y}$. If the $Y_j$'s are topological spaces and the maps
continuous maps, the direct limit set is a topological space with
the roughest topology such that all the $i_j$'s are continuous.

\smallskip

Once these standard definitions have been recalled, let us consider again an EFS $\mathcal{F}$ together with the compact set
$\Omega(\mathcal{F})= \lim_\leftarrow \mathcal{F} $ and with the following associated inverse limits $\mathcal{M}^m(\mathcal{F}) = \lim_\leftarrow
W^m(\mathcal{F},\mr) $,  $\mathcal{M}^{*m}(\mathcal{F}) =
\lim_\leftarrow W^{*m}(\mathcal{F},\mr) $. $\mathcal{M}^{*}(\mathcal{F})=\cup_{m\in \mr^+}
\mathcal{M}^{*m}(\mathcal{F})$,
$\mathcal{M}(\mathcal{F})=\cup_{m\in \mr^+}
\mathcal{M}^{m}(\mathcal{F})$.

The next proposition is a consequence of standard facts 
on the co-homology of topological projective limits.

\begin{prop}
\label{KHstandard}
$\lim_\rightarrow H^d(\mathcal{F},\ma)=H^d(\Omega(\mathcal{F}),\ma)$.
\end{prop}

Since $\mathcal{M}(\mathcal{F})=
\lim_\leftarrow  Z_d^{\geq 0}(\mathcal{F},\mr),$
there is  a natural pairing
$$\mathcal{M}(\mathcal{F})\times
H^d(\Omega (\mathcal{F}),\mr)\to \mr
\ \ \ \ \ \ <\mu|h>=h(\mu)\ $$
 defined as follows :

For $h=(h_j, \dots,h_s,\dots)$ with  $h_{s+1} = f_s^*(h_s)$, and $\mu =
(\mu_0,\dots, \mu_s,\dots)$, then $<\mu|h>=<\mu_j|h_j>$. 

The following result is a direct consequence of lemma \ref{integral}.

\begin{cor}
\label{progap}
Let $\mu= (\mu_1, \dots, \mu_n, \dots)$ be  in $\M(\F)$,
$c =(c_0, \dots , c_n, \dots)$ be a class in $H^d(\Omega(\F), \mz)$
and let $s$ be big enough so that $<\mu, c> = <\mu_s, c_s>$ then :
$$<\mu, c>\, =\, m_1\mu_{s,1}\,+\dots\, +m_{p(s)}\mu_{s,p(s)},$$
where $B_s$ has $p(s)$ $d$-regions, $\mu_{s, i}$ is the weights of the
$i^{th}$ $\,d$-region of $B_s$ and the $m_i$'s are integers. Moreover,
all such linear combinations with integer coefficients arise in this
way.
\end{cor}

\noindent In the next section we will show that $\Omega 
(\mathcal{F})$ is actually a space of tilings and that every continuous hull
$\Omega_T$ can be obtained in this way. We will also see how
$\omega_{\mathcal{F}}= \lim_\leftarrow {\rm Par}\ (\mathcal{F})$
makes $\Omega (\mathcal{F})$ a dynamical system
supporting an action of $\RM^d$, which is conjugate to the usual
action of $\RM^d$ on tiling spaces. In section \ref{invariantmes}
we will identify $\mathcal{M}(\mathcal{F})$ with the set of
invariant measures on $\Omega (\mathcal{F})$. In section
\ref{sec-cstar} we will study the direct limit of the $K$-theory of
the $B_i$. The application to the gap-labelling will follow.

\subsection{From EFS to tilings}
\label{ssec-JM}
 We associate with an EFS ${\mathcal{F}}$, the inverse limit $\omega_{\mathcal{F}}= \lim_\leftarrow {\rm Par}\
(\mathcal{F}) $.

\begin{prop}
\label{paraction}
The set $\omega_{\mathcal{F}}$ is naturally
isomorphic to $\mr^d$ and acts on  $\Omega(\mathcal{F})$.
\end{prop}

{\it Proof.} Since, for any BOF-$d$- manifold $B_i$, there exists a natural isomorphism  
$\rho_{B_i}: \mr^d \to Par(B_i)$ and for any BOF-submersion $f_i: B_{i+1}\to B_i$, the induced map $df:Par(B_i)\to
Par(B_{i+1})$ is  such that $(\rho_{B_{i+1}})^{-1}\circ df \circ \rho_{B_{i} }= id,$ it follows that 
$\lim_\leftarrow {\rm Par}$ is isomorphic to the inverse limit  $\lim_\leftarrow { \mr^d}$ associated with the inverse sequence $Id:\mr^d\to \mr^d$. Thus $\omega_{\mathcal{F}}$ is isomorphic to $\mr^d$.  

\noindent Let us show now that $\omega_{\mathcal{F}}$ acts on $\Omega({\mathcal{F}})$. From Remark \ref{rayon}, we know that there exists $r>0$ such that for each point $x= (x_0, \dots, x_i, \dots)$ in  $\Omega({\mathcal{F}})$
there exists a sequence of normal neighborhoods $U_i$ with radius
$r$ around each point $x_i$ in $B_i$ such that for each $ i >0$, $f(U_i)$ is one single sheet $D_{i-1}$ of  $U_{i-1}$ and thus $f(D_i) = D_{i-1}$. This gives a meaning to the
notion of ``small'' translation of the point $(x_0, \dots, x_i,
\dots )$. Finally any vector $v$ in $\omega_{\mathcal{F}}$ can be decomposed in a sum $v=\sum_{l=1}^{l=m}u_l$ where $\Vert u_l\Vert <r/2$ for $l=1, \dots, m$. For a point $x$ in 
$\Omega(\mathcal{F})$, we define $x+v = (\dots (x+u_1)+u_2)+\dots )+u_m)$. It is plain to check that this definition is independent of the decomposition.

\medskip

Actually  the dynamical system $(\Omega(\mathcal F), \omega_{\mathcal F})
$ is a  tiling space. Each region of the first BOF-d-manifold $B_0$ is a well oriented $d$-rectangle. Let $X$ denote the set of prototiles made with all these well oriented $d$-rectangles.
With any point $x= (x_0, \dots, x_n,\dots)$ in
$\Omega(\F)$ we associate a tiling $h(x)$ in $T(X)$ made with  prototiles in $X$ as follows:

\begin{itemize}
\item Since the injectivity radius of the $B_i$'s goes to infinity with $i$, there exists, for $i\geq 0$ a sequence of normal neighborhoods $U_i$ of $x_i$ with radius $r_i$ where $r_i$ is an increasing sequence going to infinity with $i$. 

\item From this sequence of normal neighborhoods we extract a sequence of (pre-images by  chart maps of) sheets $D_i\subset U_i$ such that $D_i\subset f_i(D_{i+1})$. These sheets are well oriented $d$-cubes centered at $x_i$.

\item Consider the translated copies $D'_i = D_i - x_i$. The $d$-cubes $D'_i$ are centered at $0$ in $\RM^d$, they have an increasing radius going to infinity with $i$. Furthermore all these $d$-cubes are tiled with the prototiles in $X$ and, for each $i>0$, the tiling of $D'_{i}$ coincides with the tiling of $D'_{i-1}$ in $D'_{i-1}$. 

\item the limit of this process defines a single tiling $h(x)$ in $T(X)$.
\end{itemize}

\noindent The proofs of the following properties are plain.

\begin{prop}
\begin{itemize}

 \item the map $h:\Omega(\F)\to T(X)$ is continuous, in\-jec\-ti\-ve 
 and conjugates the dynamical
 systems $(\Omega(\F), \omega_\F)$ and $h(\Omega(\F), \omega)$
 where  $\omega$ stands for the restriction of the standard $\mr^d$-action on $T(X)$;

 \item by construction all tilings in $h(\Omega(\F))$ satisfy 
 the finite pattern condition.

\end{itemize}
\end{prop}

\noindent As $\Omega(\F)$ is compact, it is actually the union of continuous
hulls of tiling orbits.

\medskip

So far, we have not been concerned with the minimality of the
dynamical system $(\Omega(\F), \omega_\F)$; we give now a simple
criterion that insures minimality. An EFS $\F$ verifies the {\it
repetitivity condition} if for each $n\geq 0$ and each flat
(pre-image by  chart map of a) sheet
 $D_i$ in $B_i$, there exists $p>0$ such that each
$d$-region of $B_{i+p}$ covers $D_i$ under the composition $f_i\circ
f_{i+1}\circ\dots \circ f_{i+p-1}$.

\begin{prop} Let $\F$ be an EFS that satisfies the repetitivity condition 
then the dynamical system $(\Omega(\F), \omega_\F)$ is minimal.
\end{prop}

{\it Proof.}
The proof is exactly the same as the one for tilings (proposition 1.4).

\medskip

So far we have shown how to associate with an EFS that satisfies the
repetitivity condition, a compact space of perfect  tilings made with rectangular prototiles, and  the continuous hull of  such a
tiling. On the other hand, we have discussed in section 2, the
correspondence between perfect tilings and minimal tilable
lamination. In order to have a complete scheme of these
correspondences it remains to us to be able to associate with any
perfect tiling  made  a with finite number of well oriented rectangular prototiles, an EFS that satisfies the repetitivity condition.

\subsection{From tilings to EFS}
\label{ssec-titoefs}

\noindent Consider a   perfect tiling $T$ made with the set of
prototiles $X =\{p_1, \dots, p_n\}$ where the $p_i$'s are well oriented $d$-rectangles. We are going to associate with
the hull $\Omega_T$ a  first BOF-$d$-manifold and then an EFS that satisfies
the repetitivity condition.

\noindent The First BOF-$d$-manifold $B_0$ is obtained by identifying two points in two different prototiles $x_1$ in $p_1$ and $x_2$ in $p_2$ if there exist a tile $t_1$ with $t$-type $p_1$ and a tile $t_2$ with $t$-type $p_2$ such that the translated copy of $x_1$ in $t_1$ coincides with the translated copy of $x_2$ in $t_2$. There exists a natural projection $\pi_0:\Omega_T\to B_0$
defined as follows: Consider a tiling $T'$ in $\Omega_T$, consider the
prototile where the point $0$ is and its location in the
prototile. This defines a single point in $B_0$. The  projection $\pi_0$ also defines a box decomposition ${\cal B}^{(0)}$ (see Definition \ref{boite})  of the hull equipped with its tilable lamination structure $(\Omega_T, {\cal L})$, where the sets $ \pi_0^{-1} (int (p_i))$, for $i=1, \dots, n$ are the $n$ boxes of the box decomposition.

\noindent The BOF-$d$-manifold $B_0$ is a very poor approximation of the continuous 
hull $\Omega_T$.  In order to define better and
better approximants of $\Omega_T$, we construct a sequence of box decompositions  
${\cal B}^{(i)}$, $\,{i\geq 0}$  such that for each $n\geq 0$, ${\cal B}^{(i+1)}$ is zoomed out of ${\cal B}^{(i)}$ and forces its border (see Corollary \ref{corect}) and associate with this sequence a nested sequence of tilings $(T^{(i)})_{i\geq 0}$ where $T^{(0)} = T$ (see Definition \ref {nested}). For $i\geq 1$, we construct the BOF-$d$-manifold $B_i$ as follows:

\begin{itemize}

\item First we construct the BOF-$d$-manifold $B'_i$ associated with $T^{(i)}$ in the same way we constructed $B_0$ associated with $T_0$ and a projection $\pi'_i:\Omega_T\to B'_i$. 

\item Since the tiles of $T^{(i)}$ are not well oriented $d$-rectangles but blocks and   the BOF-$d$-manifold $B'_i$ is not necessarily cellular. 

\item Each region of the BOF-$d$-manifold $B'_i$ is tiled with the prototiles in $X$. The BOF-$d$-manifold $B_i$ coincides $B'_i$ but tiled with the prototiles in $X$. The $d$-regions of $B_i$ are prototiles in $X$. 

\end{itemize}

It is clear that for $i\geq 0$, $B_i$ is a cellular BOF-$d$-manifold whose $d$-regions are well oriented $d$-rectangles, the map $\pi'_i:\Omega_T\to B'_i$ induces naturally a map $\pi_i:\Omega_T\to B_i$ and there
exists  a canonical map $f_i: B_{i+1}\to B_i$ that is defined by
$$f_i(x) = \pi_i( \pi_{i+1}^{-1}(x)).$$

Notice that this last definition makes sense since the vertical of a point of $(\Omega_T, {\cal L})$ in a box of ${\cal B}^{(i+1)}$ is included in its vertical  in a box of ${\cal B}^{(i)}$ (see Definition \ref{zoomout}). It is clear that, for each $i\geq 0$, the map $f_i:B_{i+1}\to B_i$ is a BOF-submersion and satisfies

$$f_i\circ \pi_{i+1}\,=\, \pi_i.$$

Furthermore, the sequence is an EFS since the flattening condition is a direct consequence of the fact that each box decomposition ${\cal B}^i$, for $i\geq 1$ forces its border (see for instance \cite{KP00}). Let us denote by $\F$ the EFS we have gotten, $\Omega(\F)
=\lim_\leftarrow \F$ and $\pi:\Omega_T\to \Omega(\F)$ the map
defined by $\pi(x)= (\pi_0(x), \pi_1(x),\dots, \pi_n(x),\dots)$.

\begin{prop}The map $\pi:\Omega_T\to \Omega(\F)$ is a 
conjugacy between $(\Omega_T, \omega)$
and $(\Omega(\F), \omega_\F)$.
\end{prop}

{\it Proof.} The proof is straightforward.

\medskip

\section{Invariant measures}
\label{invariantmes}
 So far we have described a same minimal dynamical system in 3 
different ways:

\begin{itemize}

\item either as the continuous hull of a perfect tiling;

\item or as a minimal tilable lamination;

\item or as a minimal expanding flattening sequence.

\end{itemize}

\noindent The interplay between these 3 points of view  is going to  provide 
a combinatorial description of the invariant measures of this dynamical 
system. Before dealing with this very precise context, let us recall first 
some basics about transverse invariant measures for oriented laminations in the large 
(we refer to\cite{Ghys} for a more complete description).

\noindent Consider a lamination $(M, \L)$ with $d$-dimensional leaves and 
fix an atlas in $\L$ made with the charts 
$h_i:U_i\to V_i\times T_i$ where $V_i$ is an open set in $\RM^d$
and $T_i$ is some topological space. Recall that  the 
transition maps $h_{i,j}\, =\, h_j\circ h_i^{-1}$ read on their
domains of definitions: 
$$h_{i, j}(x, t)\,=\, (f_{i,j}(x, t),\,\gamma_{i, j}(t)),$$ 
where $f_{i,j}$ and $\gamma_{i, j}$ are
continuous in the $t$ variable and $f_{i, j}$ is
smooth in the $x$ variable.

\begin{defi}
Let $(M, \L)$ be a lamination. 
A {\it finite transverse invariant measure} on $(M, \L)$ is the data of a 
finite positive measure on each set $T_i$ in such a way that if $B$ is a 
Borelian set in some $T_i$ which is contained in the definition set of the 
transition map $\gamma_{ij}$ then:
$$\mu_i(B)\, =\, \mu_j(\gamma_{ij}(B)).$$
\end{defi}

\noindent It is clear that the data of a transverse invariant measure for a 
given atlas provides another invariant measure for any equivalent atlas and 
thus it makes sense to consider a transverse invariant measure $\mu^t$ of a 
lamination. The fact that the leaves of a lamination carry a structure of 
$d$-dimensional manifold allows us to consider differential forms on the 
lamination $(M, \L)$. A $k$-differential form on $(M, \L)$  is the data of 
$k$-differential forms on the open sets $V_i$ that are mapped one onto the 
other by the differential of the transition maps $f_{ij}$. We denote by 
$A^k(M, \L)$ the set of $k$-differential forms on $(M, \L)$.

\begin{defi}
A {\it foliated cycle} is a linear form from $A^d(M, \L)$ to $\mr$ which is 
positive on positive forms and vanishes on exact forms.
\end{defi}

\noindent There exists a simple way to associated with a transverse invariant 
measure  a foliated cycle. Consider a $d$-differential form $\omega$ in 
$A^d(M, \L)$ and assume for the time being, that the support of $\omega$ is 
included in one of the $U_i$'s. In this case, the form can be seen as a form 
on $V_i\times T_i$. By integrating $\omega$ on the slices $V_i\times \{t\}$ 
we get a real valued map on $T_i$ that we can integrate against the transverse 
measure $\mu_i$ to get a real number $\C_{\mu^t}(\omega)$. When the support 
of $\omega$ is not in one of the $U_i$'s, we choose 
 a partition of the unity $\{\phi_i\}_i$ associated with the cover of $M$ by 
the open sets $U_i$ and define:
$$\C_{\mu^t}(\omega)\,=\, \sum_i \C_{\mu^t}(\phi_i\omega).$$
It is clear that we have defined this way a linear form $ \C_{\mu^t} : 
A^d(M, \L)\to \mr$ which does not depend on the choice of the atlas in $\L$ 
and of the partition of the unity. It is also easy to check that this linear 
form is positive for positive forms. The fact that $\C_{\mu^t}$ vanishes on 
closed form is a simple consequence of the invariance property of the 
transverse measure. The foliated cycle $\C_{\mu^t}$ is called {\it the 
Ruelle-Sullivan current} associated with the transverse invariant measure 
$\mu^t$. 
It turns out that the existence of a foliated cycle implies the existence of 
a transverse invariant measure (see \cite{Sullivan}) and thus both points of 
view: transverse invariant measure and foliated cycle are equivalent.

Let us consider now the more particular case of a tilable lamination $(M, 
\mathcal{L})$. Recall that on   a tilable lamination it is possible to
define a parallel transport along each leaf of the lamination which
gives a sense to the notion of translation in the lamination by a
vector $u$ in $\RM^d$, defining this way a dynamical system $(M,
\omega_{\mathcal{L}})$. Let $\mu $ be a measure of finite mass on $M$ 
invariant under the $\RM^d$-action. This invariant measure defines a 
transverse invariant measure of the lamination as follows. For any Borelian 
subset of a transverse set $T_i$:
$$\mu_i(B)\, =\, \lim_{r\to 0^+}\frac{1}{\lambda_d(V_i)}\mu (h_i^{-1}( V_i\times B)),$$ 
where $\lambda_d$ stands for the Lebesgue measure in $\RM^d$.  Conversally,
consider a transverse invariant measure $\mu^t$ of the tilable lamination  
$(M, \mathcal{L})$. Let $f: M\to \mr$ be  a continuous function  and assume 
for the time being that the support of $f$ is included in one of the 
$U_i$'s. In this case, the map $f\circ h_i^{-1}$ is defined on $V_i\times 
T_i$. By integrating $f\circ h_i^{-1}$ on the sheets $V_i\times \{t\}$ 
against the Lebesgue measure $\lambda_d$ of $\RM^d$, we get a real valued 
map on $T_i$ that we can integrate again the transverse measure $\mu_i$ to 
get a real number $\int f d\mu$. When the support of $f$ is not in one of 
the $U_i$'s, we choose 
 a partition of the unity $\{\phi_i\}_i$ associated with the cover of $M$ by 
the open sets $U_i$ and define:
$$\int f d\mu\,=\, \sum_i \int f\phi_i d\mu.$$
It is clear that we have defined this way a finite  measure on $M$ which 
does not depend on the choice of the atlas in $\L$ and of the partition of 
the unity and is invariant under the $\RM^d$-action.
It is also plain that the existence of a finite measure on $M$ invariant for the 
$\RM^d$-action is in correspondence with a finite  measure on a transversal 
$\Gamma$ invariant under the action of the holonomy groupoid.

Thus, for a tilable lamination  $(M, \mathcal{L})$ the following 4 points of 
view are equivalent:

\begin{itemize}

\item A finite transverse invariant measure;

\item a foliated cycle;

\item a finite  measure on $M$ invariant for the $\RM^d$-action.

\item a finite  measure on a transversal $\Gamma$ invariant for the 

holonomy groupoid action.

\end{itemize}

   Let us now develop what  follows for tilings from the above discussion. 
To fix the notations, we start with a perfect tiling $T$ made with prototiles which are well oriented $d$-rectangles  $X =\{p_1, \dots , p_q\}$. We associate with this tiling an EFS $\F$ constructed with a sequence of BOF-manifolds $B_n$, $n\geq 0$, and 
BOF-manifolds submersions: $f_n:B_{n+1}\to B_n$ that satisfy the flattening 
condition.

Let us denote by $\M(\Omega_T, \omega)$ the set of finite measures
on $\Omega_T$ that are invariant under the $\RM^d$-action. We
denote by $\M^m(\Omega_T, \omega)$ the set of finite measures on
$\Omega_T$ with total mass $m$.

Choose a set $Y$ made of one point in the interior of each
prototiles. The transversal $\Omega_{T, Y}$ is a Cantor set on
which acts the holonomy groupoid $\H_{T, Y}$. Let us denote by
$\M(\Omega_{T, Y},\H_{T, Y})$ the set of finite measures on
$\Omega_{T, Y}$ that are invariant under the action of the
holonomy groupoid $\H_{T, Y}$. With any finite invariant measure
$\mu$ in $\M(\Omega(T), \omega)$ can be associated a finite {\it
transverse}  measure $\mu^t$ in $\M(\Omega_{T, Y},\H_{T, Y})$ and
this map is one to one.  
Since $\Omega_{T, Y} $ is a Cantor set, it can be cover by a
partition in clopen sets with arbitrarily small diameters. Such a
partition $\P$ is {\it finer} than another partition $\P'$ if the
defining clopen sets of the first one are included in clopen sets
of the second one. Consider a sequence of partitions $\P_n$,
$n\geq 0$ of $\Omega_{T, Y} $ such that, for all $n\geq 0$ ,
$\P_{n+1}$ is finer than $\P_n$ and the diameter of the defining
clopen sets of $\P_n$ goes to zero as $n$ goes to $+\infty$.

\noindent{\bf Claim} {\it A finite measure on $\Omega_{T, Y} $ is given
by the countable data of non negative numbers associated with each defining
clopen sets of each partition $\P_n$ which satisfy the obvious
additivity relation.}

The EFS $\F$ provides us such a sequence of partitions $\P_n$ as
follows. For each $n\geq 0$, let $F_1, \dots, F_{p(n)}$ be the
regions of the BOF-$d$-manifold $B_n$. Each $F_i$ is a copy of a prototile in $X$ and consequently there is a marked point $y_i$ in $F_i$. For $n\geq 0$, $i=1, \dots , p(n)$, we consider the clopen set:
 $$\C_{n, i}\, =\,\pi_n^{-1}(y_i).$$
 For $n$ fixed and when $i$ varies from $1$ to
$n(p)$ the clopen sets $\C_{n, i}$ form a partition $\P_n$ of
$\Omega_{T, Y} $. Furthermore, for $n\geq 0$,
$\P_{n+1} $ is finer than $\P_n$ and  the diameter of the
clopen sets $\C_{n,i}$ goes to zero as $n$ goes to $+\infty$. It
follows that a finite measure on $\Omega_{T, Y} $ is  given by the
countable data of non negative weights  associated with each defining clopen
sets  $\C_{n,i}$ which satisfy the obvious additivity relation.

The relation between an invariant measure $\mu$ in
$\M(\Omega_T, \omega)$  and the associated transverse invariant measure 
$\mu^t$ in $\M(\Omega_{T, Y},\H_{T, Y})$  is given by :

\begin{prop}\label{trans}
For $n\geq 0$ and $i$ in $\{1, \dots n(p)\}$:
$$\mu^t(\C_{n, i})\,=\, {\frac{1}{\lambda_d(F_i)}}\mu (\pi_n^{-1}(F_i)).$$
\end{prop}

In the sequel, we are going to characterize in a combinatorial way these 
invariant measures.

\noindent Consider the map $\tau_n: \M(\Omega_T, \omega)) \to C_d(B_n, \mr)$ defined 
by
$$\tau_n(\mu) = (\frac{\mu(\pi_n^{-1}(F_1))}{\lambda_d(F_i)}, \dots, 
\frac{\mu(\pi_n^{-1}(F_{p(n)}))}{\lambda_d(F_{p(n)})}),$$
 where the regions
$F_1, \dots, F_{p(n)}$ are now ordered and equipped with the natural 
orientation that
allows to identify $C_d(B_n, \mr)$ with $\mr^{p(n)}$.

\begin{remark} \label{coef} Notice that the coordinates of $\tau_n(\mu)$ are 
the transverse measures associated with $\mu$ of clopen sets $\C_{n,io}$ for 
some  $n\geq 0$ and some $i$ in $\{1, \dots n(p)\}$.
\end{remark}

\begin{prop} \label{meas}
For any $n\geq 0$, the map $\tau_n$ satisfies the following properties:

\begin{itemize}

\item [(i)]  $\tau_n(\M(\Omega_T, \omega)) \subset W^{\star }(B_n)$;

\item [(ii)] $f_{n,*}\circ \tau_n = \tau_{n+1}$.

\end{itemize}\end{prop}

{\it Proof.}

\noindent (i) Choose an invariant measure $\mu$ in $\M(\Omega_T, \omega)$.
The invariance of $\mu$ implies that on each edge of $B_n$ the sum
of the transverse measures associated with the regions on one side of the edge
is equal to the sum of the transverse measures of regions on the
other side of the edge. These are exactly the switching rules (or
Kirchoff-like laws) that define $W(B_n)$. The fact that the measure is
invariant implies that each region has a strictly positive weight.
Thus $\tau_n(\mu)$ is in $W^{\star }(B_n)$.

\noindent (ii) Let $F'_1, \dots, F'_{p(n+1)}$ be the ordered sequence of
regions of $B_{n+1}$ equipped with the natural orientation that allows
to identify $C_d(B_{n+1}, \mr)$ with $\mr^{p(n+1)}$. To the linear
map $f_{n,*}:C_d(B_{n+1}, \mr)\to C_d(B_n, \mr)$ corresponds a
$n(p)\times n(p+1)$ matrix $A_n$ with integer non negative
coefficients. The coefficient $a_{i, j, n}$   of the $i^{th}$ line
and the $j^{th}$ column is exactly the number of pre images in
$F'_j$ of a point in $F_i$.  Thus we have the relations:
$$\frac{\mu(\pi_n^{-1}(F_i))}{\lambda_d(F_i)}\, =\, \sum_{j=1}^{j=p(n+1)}a_{i, j, n}\frac{\mu(\pi_{n+1}^{-1}(F'_j))}{\lambda_d(F'_j)},$$ 
for all $i=1,\dots, p(n)$ and
all $j=1, \dots p(n+1),$ which exactly means that $f_{n,*}\circ
\tau_n = \tau_{n+1}.$

\medskip

\noindent Let us set $\M^{\star }(\F) = \lim_\leftarrow W^{\star }$. Then we 
have the following characterization of the set of invariant measures of 
$(\Omega_T, \omega)$.

\begin{teo}

$$\M(\Omega_T, \omega)\, \cong\, \M^{\star }(\F).$$
\end{teo}

{\it Proof.} The inclusion $$\M(\Omega_T, \omega)\, \subset\,
\M^{\star }(\F)$$ is a direct consequence of Proposition
\ref{meas}.

\noindent With any element $(\beta_0, \dots, \beta_n, \dots)$ in
$\M^{\star }(\F)$ we can associate, thanks to Proposition
\ref{trans} a weight to each clopen set $\C_{n, i}$. The
relation $\beta_{n} =f_{\star n}\beta_{n+1}$ means that this
countable sequence of weights satisfies the additivity property
and then defines a measure on $\Omega_{T, Y} $. The fact all the
$\beta_n$'s are cycles, i.e satisfy the switching rules, yields
that this measure is a transverse invariant measure, i.e. an
element in  $\M(\Omega_{T, Y},\H_{T, Y})$. Since the
correspondence between $\M(\Omega_{T, Y},\H_{T, Y})$ and
$\M(\Omega_{T},\omega)$ is bijective, the equality is proved.

\medskip

\begin{cor}

\begin{itemize}

\item If the dimension of $H^d(B_n, \mr)$ is uniformly bounded 
by $N$, then  for all $m>0$, $\M^m(\Omega_T, \omega)$ contains at most $N$
ergodic measures;

\item if furthermore the coefficients of all
the matrices  $f_{\star n}$ are uniformly bounded then
for all $m>0$, $\M^m(\Omega_T, \omega)$ is reduced to a single point i.e;
the dynamical system $(\Omega_T, \omega)$ is uniquely ergodic.

\end{itemize}
\end{cor}

{\it Proof.} The proof is standard and can be found in \cite {G-M} in a
quite similar situation in the particular case when $d=1$. To prove the first statement we may assume that the dimension of the $H^d(B_n, \mr)$'s is constant and equal to $N$. The set $\M^m(\Omega_T, \omega)$ is a convex set and its extremal points coincides with the
set of ergodic measures. Since $\M^m(\Omega_T, \omega)\, =\, \M^{\star m}(\F)$, the convex set $\M^m(\Omega_T, \omega)$ is the intersection
of the convex nested sets:
$$\M(\Omega_T, \omega)\, =\, \cap_{n\geq 0}W_n$$ where
$$W_n\, =\, f_{\star 1}\circ \dots\circ f_{\star n-1}W^\star(B_n).$$ Since
each convex cone $W_n$ possesses at most $N$ extremal lines, the limit
set $\M^m(\Omega_T, \omega)$ possesses also at most $N$ extremal
points and thus at most $N$ ergodic measures.

\noindent In order to prove  the second statement, we have
to show that $\M(\Omega_T, \omega)$ is one dimensional. Consider
two points $x$ and $y$ in the positive cone of $\mr^N$.
Let $T$ be the largest line segment containing $x$ and $y$ and
contained in the positive cone of $\mr^N$. We recall that the
hyperbolic distance between $x$ and $y$ is given by:
$$Hyp (x, y) \, =\, -\ln{\frac{(m+l)(m+r)}{l.r}},$$ where $m$ is the
length of the line segment $[x, y]$ and $l$ and $r$ are the length of
the connected components of $T\setminus [x, y]$. Positive matrices
contract the hyperbolic distance in the positive cone of $\mr^N$. Since
the matrices corresponding to the maps $f_{\star n}$ are uniformly
bounded in sizes and entries, this contraction is uniform. Because of
this uniform contraction the set $\M(\Omega_T, \omega)$ is one
dimensional.

\begin{remark} If it is easy to construct perfect tilings in dimension $d 
=1$ which are not uniquely ergodic (see for instance \cite{G-M}), this 
question remains unclear in dimension $d\geq 2$.
\end{remark}

Let us finish this section by giving a more specific interpretation of
the Ruelle-Sullivan current associated with a transverse invariant
measure in the case of tilings. Since the Ruelle-Sullivan current
vanishes on exact $d$-differential forms, it acts on the de-Rham
cohomology group $H^d_{DR}(\Omega (\F))$. By the same standard
arguments as the one developed in Section \ref{ssec-EFS} we have:
$$ H^d_{DR}(\Omega (\F))\,=\,\lim_\rightarrow H^d_{DR}(\mathcal{F}).$$
In other words, every co-homology class $[\omega]$ in $H^d_{DR}(\Omega
(\F))$ is the direct limit of co-homology classes $[\omega_n]$ in
$H^d_{DR}(\Omega (\F))$. For $n$ big enough, we have the relation:
$$\C_{\mu^t}([\omega])\, =\,\sum_{i=1}^{p(n)}\frac{\mu(\pi_n^{-1}(F_i))}{\lambda_d(F_i)}\int_{F_i}\omega_n.$$
 Let $I_n$ be the standard isomorphism $I_n: H^d_{DR}(B_n)\to H^d(B_n, \ma)$ (where $\ma =\mr$ or $\mc$ according
to the fact that the forms have coefficients in $\mr$ or $\mc$)
defined by:
$$<I_n([\omega)], c>\, =\, \int_c\omega,$$ for every cycle $c$ in
$H_d(B_n, \ma)$. 

\noindent As we already suggested it will be important for the
proof of the gap-labelling theorem to consider the integral
co-homology classes. The set of integral classes on $B_n$ is defined
by:
$$H^{d, int}_{DR}(B_n)\, =\,I_n^{-1}H^d(B_n, \mz),$$ it is the set of
classes that take integer values on integers cycles (i.e; cycles in
$H_d(B_n, \mz)$.  We define:
$$H^{d, int}_{DR}(\Omega (\F))\, =\, \lim_\rightarrow H^{d, int}_{DR}(\F).$$ 

From Corollary \ref{progap} and Remark \ref{coef} we get:

\begin{prop}\label{quasigap}
$$\C_{\mu^t}(H^{d, int}_{DR}(\Omega (\F))) = \int_{\Omega_{T, Y}} d\mu^t \; \Cc (\Omega_{T, Y}, \ZM)$$
 where $\Cc (\Gamma_T, \ZM)$ is the set of continuous functions on $\Omega_{T, Y}$ with integer values.
\end{prop}

\medskip

\section{$C^*$-algebras, $K$-theory and Gap-Labelling}
\label{sec-cstar}

\noindent In the previous section we have seen how the description of 
 the continuous hull $\Omega_T$  of a perfect tiling of $\RM^d$ in terms of expanding 
flattening sequences has been powerful in order to describe 
the ergodic properties of the $\RM^d$-action on $\Omega_T$.
In this section we are going to use the same tools to analyze  the $K$-theory of 
$\Omega_T$.

\medskip

Let us remind some very elementary facts of classical ``topological''
$K$-theory (see for instance \cite{At}).

\noindent For any Abelian semigroup $(A,+)$ with zero $0$, there is a
canonical way to associate with $A$ an Abelian group $(K(A),+)$
(also called de {\it Grothendieck} group of $A$) which satisfies a
natural universal property. A simple way to construct
$K(A)$ is as follows. Consider the product semigroup $A\times A$
and let $\Delta$ be its diagonal. The cosets $[(a,b)]=(a,b)+
\Delta$ make a partition of $A\times A$ and we can define on the
cosets set $K(A)=A\times A/\Delta$ the operation
$[(a,b)]+[(c,d)]=[(a+c,b+d)]$. Note that $[(a,b)]+[(b,a)]=[(0,0)]$,
hence $K(A)$ an Abelian group. $$\alpha : A\to K(A) \ \ \
\ \ \alpha (a)=[a]=[(a,0)]$$ is a semigroup homomorphism which
satisfies the following universal property: for any group $G$ and
semigroup homomorphism $\gamma : A\to G$ there exists a unique
homomorphism $\chi : K(A)\to G$ such that $\gamma = \chi \alpha$.  If
$A$ satisfies the {\it cancelation rule} (if $a+b=c+b$ then $a=c$),
then $\alpha$ is injective. Let us denote by $K^+(A)$ the image of
$\alpha$. This is a ``positive cone'' in $K(A)$ (with respect to the
$\mz$-module structure of $K(A)$) which generates the whole
$K(A)$. Every $[(a,b)]$ is of the form $[(a,b)]=[a]-[b]$. If for every
$a,b\in A$, $a+b=0$ if and only if $(a,b)=(0,0)$, then $K^+(A)\cap
-K^+(A)= \{0\}$, hence the relation $[a]-[b]\geq [a']-[b']$ iff
$[a]-[b]-[a']+[b']\in K^+(A)$ makes $K(A)$ an ordered group.

\smallskip

\noindent Let $X$ be a {\it compact} topological space.  Apply the above
construction to the semigroup of isomorphism classes of complex vector
bundles on $X$ with the operation given by the direct sum
$\oplus$. The class of the unique rank $0$ vector bundle is the $0$ of
this semigroup.  The resulting Abelian group is denoted $K^0(X)$. So
each element of $K^0(X)$ is of the form $[E]-[F]$, where $E$ and $F$
are (classes of) vector bundles on $X$; $[E]=[F]$ iff there exists $G$
such that $E\oplus G = F\oplus G$. $K^+(X)=(\{[E]-[0]\}$ and it is the
positive cone of an actual order on $K(X)$. On the other hand the semigroup
does not satisfies the cancelation rule.

Since $X$ is compact,  for every $F$ there exist $n\in \mn$ and $G$ such that $[F]\oplus [G]=\epsilon^n$, where $\epsilon^n$ denotes the (class of) trivial
vector bundle of rank $n$. Hence $[E]-[F] = [E\oplus G]- [\epsilon^n]$, so that each element $\beta$ of $K^0(X)$ is of the
form $\beta =[H]-[\epsilon^n]$, for some $n\in \mn$. If $n_0$ is
the minimum of such $n$, then $\beta \geq 0$ iff $n_0=0$.  Moreover
if $G\oplus G'=\epsilon^m$, then $E\oplus G = F\oplus G$ implies
that $E\oplus \epsilon^m = F\oplus \epsilon^m$. We can summarize
this remark by saying that $[E]=[F]$ if and only if $E$ and $F$
are {\it stably equivalent}. 

\smallskip

\noindent If $\Cc(X)$ denotes the ring of
continuous complex valued functions on $X$, we can associate with
each vector bundle $E$ on $X$ the $\Cc (X)$-module of the sections of
$E$, $\Gamma(E)$. This is in fact a functor $\Gamma$ from the category
$\mathcal{B}$ of vector bundles over $X$ to the category
$\mathcal{M}$ of $\Cc (X)$-modules. The functor $\Gamma$ induces an equivalence
between the category $\mathcal{T}$ of trivial vector bundles to
the category $\mathcal{F}$ of free $\Cc (X)$-modules of finite rank.
As $X$ is compact, the fact that for every bundle $E$ there
exists a bundle $G$ such that $E\oplus G$ is trivial, means that
the category of vector bundles over $X$ coincides with the
sub-category $Proj(\mathcal{T})$ (which a priori is smaller)
generated by the images of {\it projection operators} on trivial
bundles. $Proj(\mathcal{F})$ is defined in a similar way. Hence
$\Gamma$ establishes an equivalence between $Proj(\mathcal{T})$ to
$Proj(\mathcal{F})$ which is by definition the category of {\it
finitely-generated projective} $\Cc (X)
$-modules. The construction of
$K^0(X)$ can be rephrased in terms of projective modules instead of
vector bundles and in such a case it is denoted $K_0(\Cc (X))$.

\smallskip

\noindent The Chern classes of a bundle depends only on its equivalence class  $\beta $ in $K^0(X)$, we denote them by $c_i(\beta)\in H^{2i}(X,\mz)$. If $f: X\to Y$ is a continuous map between compact spaces, then there is a natural map $f^*_k: K^0(Y)\to K^0(X)$, such that for any $\beta \in K(Y)$, $c_i(f^*(\beta))=f^*(c_i(\beta))$, where $f^*: H^{2i}(Y,\mz)\to H^{2i}(X,\mz)$ is the natural map induced by $f$ on co-homology.

\smallskip

\noindent The group $K^1(X)$ denotes by definition the group $K^0(SX)$, where 
$SX$ is the {\it reduced suspension } of $X$ , i.e. the quotient space $X\times {\bf S}^1/ 
\{x\}\times{\bf S}^1\cup X\times \{p_1\}$ where $x$ is a marked point in $X$ and 
$p_1$ a marked point in the circle ${\bf S}^1$ ( see \cite{At}). Let $p: X\times 
{\bf S}^1\to SX$ be
the natural projection. If $f:X\to Y$ is as above, $(f\times id): X\times S^1 \to Y\times {\bf S}^1$ induces a continuous map  $Sf: SX\to SY$ which can be used
to define $f^*_k: K^1(Y)\to K^1(X)$.

\smallskip

Let us now come back to the context of tilings. Since
$K^0$ and $K^1$ have a good behavior with respect to projective
limits, we associate with an  EFS
$\mathcal{F}=\{(B_n,f_n)\}_{n\in \mn}$, the direct sequence  
$K^i(\mathcal{F})=\{(K^i(B_n),(f^*_k)_n)\}_{n\in \mn}$, $i=0,1$.  We have: 
$$K^i(\Omega_T)=\lim_\rightarrow K^i(\mathcal{F})\ .$$

\noindent The last statement can be understood also in terms of 
projection operators. Let $\psi_n:\Omega_T \mapsto B_n$ be the
canonical map defined through the projective limit. Therefore any
continuous function $f\in\Cc(\Omega_T)$ is the uniform limit of a
sequence of functions $f_n\in \Cc (B_n)$ namely $\lim_{n\uparrow \infty} \| f-f_n \circ \psi_n \| =0$, the norm being the sup-norm. In particular, any
projection $P\in \Cc (\Omega_T) \otimes M_N(\CM)$ is equivalent to a
projection in $\Cc (B_n) \otimes M_N(\CM)$ for at least one $n\in
\NM$, which means  that the $K$ groups of $\Omega_T$
is the inductive limit of the $K$-groups of the $B_n$'s.

Let us remind now some facts about the $C^*$-algebra $\Aa = \Aa_T$
which represents the deep non-commutative geometric structure of our
dynamical system $\Omega_T$. The $C^*$-algebra $\Aa$ is the crossed product
\Cs $\Cc (\Omega_T)\rtimes \RM^d$. Let $\Aa_0$ be the dense sub algebra
made of continuous functions on $\Omega_T\times \RM^d$ with compact
support. If $\mu$ is an $\RM^d$-invariant probability measure on
$\Omega_T$, there is a canonical trace $\TVm$ defined on $\Aa_0$ by
$\TVm (A) = \int_{\Omega_T} d\mu(\omega) A(\omega, 0)$.  The main
object of the {\it gap-labelling} question concerns the trace $\TV_{\mu} (P)$
of a projection $P\in\Aa$.

\medskip

\noindent The definition of $K_0$ recalled above in the classical topological
setting can be generalized to non commutative algebras. Recall
\cite{Black86} that two projections $P,\,Q \in
\Aa$ are equivalent whenever there is an element $U\in \Aa$ such
that $P=UU^{\ast}$ and $Q= U^{\ast}U$. If $P$ and $Q$ are
orthogonal to each other, the equivalence class of their direct
sum $[P\oplus Q]$ depends only upon $[P]$ and $[Q]$, leading to
the definition of the addition $[P]+[Q]=[P\oplus Q]$. To make sure
that two projections can be always be made mutually orthogonal
modulo equivalence, $\Aa$ must be replaced by $\Aa\otimes \Kk$.
Here $\Kk$ is the \Cs of compact operators on a Hilbert space with
a countable basis. It can also be defined as the smallest \Cs
containing the increasing sequence of finite dimensional matrices:
$\Kk =\lim_{\rightarrow} M_n(\CM)$, where the inclusion of $M_n$
into $M_{n+m}$ ($0<m$) is provided by
$$A\in M_n \mapsto i_{n, m}(A) \;=\;   \left[\begin{array}{cc}
A & 0\\
 0 & 0_{m}
   \end{array}\right]\in M_{n+m}
$$

\noindent The group $K_0(\Aa)$ is the group generated by 
formal differences $[P]-[Q]$
of equivalent classes of projections in $\Aa\otimes \Kk$.  Then two
equivalent projections have the same trace, and since the trace is
linear, it defines a group homomorphism $\TV_{\mu,\ast} : K_0(\Aa)
\mapsto \RM$ the image of which are called the {\em gap labels}.

\noindent Together with $K_0$, there is $K_1$ which is 
defined as the equivalence classes, under homotopy, of invertible elements in $\lim_{\rightarrow}
GL_n (\Aa)$. In the case of $K^1(\Cc (X))$ this is equivalent to the previous
definition in terms of reduced suspensions.  
A standard result due to  R. Bott \cite{Black86} claims that
$K_1(\Aa)$ is isomorphic to $K_0(\Aa\otimes \Cc_0 (\RM))$ and that
$K_1(\Aa\otimes \Cc_0 (\RM))$, which is then nothing but $K_2(\Aa) =
K_0(\Aa\otimes \Cc_0 (\RM^2))$, is actually isomorphic to $K_0(\Aa)$
(Bott's periodicity theorem). Both $K_0$ and $K_1$ are discrete
Abelian groups. They are countable whenever $\Aa$ is separable. An
important property of $K(\Aa)=K_0(\Aa)\oplus K_1(\Aa)$ is that it
defines a covariant functor which is continuous under taking inductive
limits. Namely any $^{\ast}$-isomorphism $\alpha:\Aa\mapsto \Bb$
between \Css induces a group homomorphism $\alpha_{\ast}$ defined by
$\alpha_{\ast}([P]) = [\alpha(P)]$ for $K_0$ and similarly for
$K_1$. Moreover, $K(\lim_{\rightarrow}\A_n)= \lim_{\rightarrow}
K(\Aa_n)$.

 The Thom-Connes theorem \cite{Black86,Connes94}  states that there is a group 
isomorphism
$\phi_d$ between $K_0 (\Cc (\Omega_T)\rtimes \RM^d)$ and $K_d(\Cc
(\Omega_T))$ (with $K_i = K_{i+2}$ by Bott's periodicity
theorem). Actually the trace action makes this isomorphism more explicite:

\begin{theo}
\label{the-connesthom} \cite{Connes86} 
Let $P$ be a projection in $\Aa$ and let $[P]$ its class in $K_0(\Aa)$.

\noindent (i) If $d$ is odd, let $U$ be a unitary element of 
$\Cc (\Omega_T)\otimes \Kk$ representing $\phi_d([P])$. Let $\eta $ be the $d$-form in $H^d_{DR}(\Omega (\F))$:
$$ \eta = \tr \left( (U^{-1} dU)^d \right)\,.
$$

\noindent (ii) If $d$ is even, let $Q_{\pm}$ be a pair of projections of
$\Cc (\Omega_T)\otimes \Kk$ with $\phi_d([P])= [Q_+] -[Q_-]$. Let $\eta $ be
the $d$-form  in $H^d_{DR}(\Omega (\F), \mc)$
$$ \eta \;=\;    \tr \left( (Q_+\, dQ_+\wedge dQ_+)^{d/2} \right) - \tr \left( (Q_-\, dQ_-\wedge dQ_-)^{d/2} \right)  \,.
$$

\noindent Then, in both cases, and for a suitable (see later) 
normalization constant $k_d$, $k_d\eta$ is in $H^d_{DR}(\Omega (\F), 
\mz))$ and 
$$\TVm (P) \;=\; \C_{\mu^t}([k_d\eta]).
$$
\end{theo}

Combining the previous result with Proposition \ref{quasigap}, we get:
$$\TVm (K_0(\Aa)) \, \subset\,\int_{\Omega_{T, Y}} d\mu^t \; \Cc (\Omega_{T, Y}, 
\ZM).$$

\noindent The converse inclusion being a standard result,
the  {\it gap-labelling theorem} (see theorem
\ref{gaplab} stated in the introduction) is proved.

The integral classes defined above are related to Chern classes.  The image, say
$\beta$, of $[P]$ by the Thom-Connes isomorphism can be taken either
as an element of $K^0 (B_n)$, when $d$ is even, or of $K^1(B_n)$, 
when $d$ is odd, for $n$ large enough.  It turns out (this is
contained in the proof of the Thom-Connes theorem) that:

1) when $d$ is even $k_d\,\eta$ (which actually ``lives'' on $B_n$)
represents the Chern class $c_{[d/2]}(\beta)\in H^d(B_n,\mz)$; in fact
the choice of the normalization constant $k_d$ in theorem
\ref{the-connesthom} ensures that this class is integral. So
$$\TVm(P)= <\mu_n |c_{[d/2]}(\beta) > \ .$$

2) When $d$ is odd 
 $$\TVm(P)= <S\mu_n |c_{[(d+1)/2]}(\beta) > \ , $$
where $S\mu_n$ is defined as follows. For $z\mu_n\in Z_k(B_n,\ma)$, 
$(\mu_n\times S^1)\in Z_{k+1}(B_n\times S^1$
(it is understood that $S^1$ has the usual counterclockwise
orientation).  By using the natural projection $p_n: B_n\times S^1 \to
SB_n$, $(\mu_n\times S^1)$ induces a $(k+1)$ (singular) cycle $S\mu_n$ on 
$SB_n$,
which is called the {\it suspension} of the cycle $\mu_n$.

\medskip
\eject

\end{document}